\documentclass{amsart}

\usepackage{amsthm,amsmath,amssymb,amsfonts}
\usepackage[utf8]{inputenc}
\usepackage{verbatim}
\usepackage{marvosym}
\usepackage{stackrel}
\usepackage{graphicx}
\usepackage[svgnames]{xcolor}
\definecolor{blue(munsell)}{rgb}{0.0, 0.5, 0.69}
\usepackage{lipsum}
\usepackage{svg}
\usepackage{caption}
\usepackage{enumitem}
\usepackage{epigraph}
\usepackage{ebproof}
\usepackage[cal=boondoxo]{mathalfa}
\usepackage[all,cmtip]{xy}
\usepackage{mathtools}
\usepackage{scalerel}
\usepackage{color}
\usepackage{mathrsfs}
\usepackage{subfig}
\usepackage{framed}
\usepackage{hyperref}
\hypersetup{hypertexnames = false, bookmarksdepth = 2, bookmarksopen = true, colorlinks, linkcolor = black, citecolor = blue(munsell), urlcolor = blue(munsell), pdfstartview={XYZ null null 1}, breaklinks=true}
\usepackage{cleveref}
\usepackage{doi}
\usepackage{tikz}

\usepackage{bbold}
\usepackage{ccicons}
 
\usepackage[colorinlistoftodos]{todonotes}
\usepackage{xparse}
\usepackage{scalerel,stackengine}
\stackMath
\newcommand\reallywidehat[1]{%
	\savestack{\tmpbox}{\stretchto{%
			\scaleto{%
				\scalerel*[\widthof{\ensuremath{#1}}]{\kern-.6pt\bigwedge\kern-.6pt}%
				{\rule[-\textheight/2]{1ex}{\textheight}}
			}{\textheight}%
		}{0.5ex}}%
	\stackon[1pt]{#1}{\tmpbox}%
}

\DeclareMathSymbol{\mlq}{\mathord}{operators}{``}
\DeclareMathSymbol{\mrq}{\mathord}{operators}{`'}
\DeclareMathSymbol{\mlqq}{\mathord}{operators}{"5C}
\DeclareMathSymbol{\mrqq}{\mathord}{operators}{`"}

\DeclareDocumentCommand\issue{g}{\todo[size=\footnotesize,color = green!40]{Issue\IfNoValueF{#1}{: #1}}}
\DeclareDocumentCommand\tobedone{g}{\todo[size=\footnotesize,color = yellow!50]{To do\IfNoValueF{#1}{: #1}}}
\DeclareDocumentCommand\notationissue{g}{\todo[size=\footnotesize,color = red!30]{Notation?\IfNoValueF{#1}{: #1}}}
\DeclareDocumentCommand\doubt{g}{\todo[size=\footnotesize,color = blue!10]{Doubt\IfNoValueF{#1}{: #1}}}
\DeclareDocumentCommand\observation{g}{\todo[size=\footnotesize,color = orange!10]{Observation\IfNoValueF{#1}{: #1}}}

\usepackage[all,cmtip]{xy}

\usetikzlibrary{calc}
\usetikzlibrary{matrix,arrows,arrows.meta}
\usepackage{tikz-cd}

\usepackage{calligra} 

\newtagform{fn}{(}{)\footnotemark}

    {\endMakeFramed}

    {\endMakeFramed}

\newcommand{\vcenteredinclude}[1]{\begingroup
	\setbox0=\hbox{\includegraphics[scale=0.25]{#1}}%
	\parbox{\wd0}{\box0}\endgroup}
\makeatletter
\g@addto@macro\bfseries{\boldmath}
\makeatother

\theoremstyle{plain}
\newtheorem{thm}{Theorem}[section]
\newtheorem*{thm*}{Theorem}
\newtheorem{prop}[thm]{Proposition}
\newtheorem{lem}[thm]{Lemma}

\theoremstyle{definition}
\newtheorem{defn}[thm]{Definition}

\theoremstyle{remark}
\newtheorem{rem}[thm]{Remark}
\newtheorem{exa}[thm]{Example}
\newtheorem{notat}[thm]{Notation}

\newtheorem{quest}[thm]{Question}

\newcommand\Ind{\operatorname{Ind}}

\newcommand\Mod{\operatorname{Mod}}
\newcommand\fpmod{\operatorname{mod}}

\newcommand\Lan{\operatorname{lan}}

\newcommand\Set{\operatorname{\bf Set}}

\newcommand\Rex{\operatorname{Rex}}

\newcommand\op{\circ}

\newcommand\Sh{\operatorname{Sh}}

\newcommand\colim{\operatorname{colim}}

\newcommand\ca{\mathcal {A}}
\newcommand\cb{\mathcal {B}}
\newcommand\cc{\mathcal {C}}
\newcommand\cd{\mathcal {D}}
\newcommand\cf{\mathcal {F}}

\newcommand\cg{\mathcal {G}}

\newcommand\ce{\mathcal {E}}

\newcommand\ck{\mathcal {K}}

\newcommand\crr{\mathcal {R}}

\newcommand{\AAA}{\mathfrak{a}}
\newcommand{\BBB}{\mathfrak{b}}

\newcommand{\CCC}{\mathfrak{c}}
\newcommand{\DDD}{\mathfrak{d}}

\newcommand{\Cocont}{\ensuremath{\mathsf{Cocont}}}
\newcommand{\Cont}{\ensuremath{\mathsf{Cont}}}
\newcommand{\Lex}{\ensuremath{\mathsf{Lex}}}

\newcommand{\Cat}{\ensuremath{\mathsf{Cat}}}
\newcommand{\Img}{\ensuremath{\mathsf{Img}}}

\newcommand{\Kern}{\ensuremath{\mathsf{Ker}}}
\newcommand{\Cokern}{\ensuremath{\mathsf{Coker}}}
\newcommand{\QInj}{\ensuremath{\mathsf{QInj}}}
\newcommand{\Presh}{\ensuremath{\mathsf{Presh}}}
\newcommand{\Acc}{\ensuremath{\mathsf{Acc}}}
\newcommand{\Cnt}{\ensuremath{\mathsf{Cnt}}}

\newcommand\Pres{\ensuremath{\mathsf{Pres}}}
\newcommand{\site}{\ensuremath{\mathsf{Site}}}
\newcommand{\siteflat}{\ensuremath{\mathsf{Site}_{\mathsf{flat}}}}
\newcommand{\grothflat}{\ensuremath{\mathsf{Grt}_{\flat}}}

\newcommand{\groth}{\ensuremath{\mathsf{Grt}}}

\newcommand{\lfp}{\ensuremath{\mathsf{lfp}}}
\newcommand{\coh}{\ensuremath{\mathsf{coh}}}
\newcommand{\Qcoh}{\ensuremath{\mathsf{Qcoh}}}
\newcommand{\Sch}{\ensuremath{\mathsf{Sch}}}
\newcommand{\qcqs}{\ensuremath{\mathsf{qcqs}}}
\newcommand{\fflat}{\ensuremath{\mathsf{flat}}}
\newcommand{\Ringspaces}{\ensuremath{\mathsf{RingSpace}}}
\newcommand{\Spec}{\ensuremath{\mathsf{Sp}}}

\newcommand{\obj}{\ensuremath{\mathbb{O}}}

\newcommand{\ra}{\rightarrow}

\mathchardef\mhyphen="2D

\title{Exponentiable Grothendieck categories in flat algebraic geometry}
\author{Ivan \textsc{Di Liberti}} 

\author{Julia \textsc{Ramos Gonz\'alez}} 

\address{
	Ivan \textsc{Di Liberti}:\newline
	Institute of Mathematics\newline
	Czech Academy of Sciences\newline
	\v{Z}itn\'{a} 25, Prague, Czech Republic\newline
	\href{mailto:diliberti.math@gmail.com}{\sf diliberti.math@gmail.com}\newline
}

\address{
	Julia \textsc{Ramos Gonz\'alez}: \newline
	Institut de recherche en math\'ematique et physique \newline
	Universit\'e Catholique de Louvain \newline
	Cyclotron 2/L7.01.02, 1348 Louvain-la-Neuve, Belgium \newline
	\href{mailto:julia.ramos@uclouvain.be}{\sf julia.ramos@uclouvain.be}
}
\thanks{\doi{10.1016/j.jalgebra.2022.03.040}} 
\thanks{\footnoterule
	\vcenteredinclude{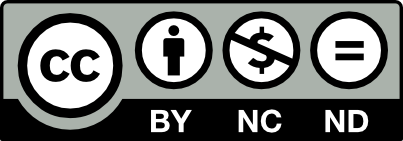} \copyright 2022. This manuscript version is made available under the CC-BY-NC-ND 4.0 license \url{https://creativecommons.org/licenses/by-nc-nd/4.0/}}
\begin{document}
\tikzcdset{arrow style=tikz, diagrams={>=to}}

\begin{abstract}  
We introduce and describe the $2$-category $\grothflat$ of Grothendieck categories and flat morphisms between them. First, we show that the tensor product of locally presentable linear categories $\boxtimes$ restricts nicely to $\grothflat$. Then, we characterize exponentiable objects with respect to $\boxtimes$: these are the continuous Grothendieck categories. In particular, locally finitely presentable Grothendieck categories are exponentiable. Consequently, we have that, for a quasi-compact quasi-separated scheme $X$, the category of quasi-coherent sheaves $\mathsf{Qcoh}(X)$ is exponentiable. Finally, we provide a family of examples and concrete computations of exponentials. 
\end{abstract}
\maketitle

\tableofcontents

\section{Introduction}

\subsection{General setting}  
From the perspective of noncommutative algebraic geometry, Grothendieck abelian categories play the role of (models of) (possibly) noncommutative schemes (see for example \cite{artinzhang94}, \cite{staffordvandenbergh01}, \cite{kontsevichrosenberg00} among many others). This intuition is motivated by the Gabriel–Rosenberg reconstruction theorem, which shows that a quasi-separated scheme can be reconstructed, up to isomorphism of schemes, solely from the abelian category of quasicoherent sheaves on the scheme, which is a Grothendieck category \cite[Tag 077P]{stacksproject}. The theorem was initially proved for noetherian schemes by Gabriel \cite[Ch. VI, \S3]{gabriel62} and generalized to quasi-separated schemes by Rosenberg \cite{rosenberg04} (see also \cite{brandenburg18}). In addition, the Gabriel-Popescu theorem \cite{gabrielpopescu64} allows to interpret Grothendieck categories as a linear-version (that is, an $\mathsf{Ab}$-version or $\Mod(\mathbb{Z})$-version) of Grothendieck topoi \cite{lowen04}, perspective that emphasizes their geometric nature.

A natural subsequent step is to try to determine which is the correct notion of morphism between noncommutative schemes, that is, which is a suitable choice of morphisms between Grothendieck categories allowing to replicate and generalize to a noncommutative setting the methods and intuition from classical algebraic geometry. Different $2$-categories of Grothendieck categories are considered in the literature with algebraic-geometric purposes. We revise here some of them:
\begin{itemize}
	\item[$\mathsf{Grt}$] is the $2$-category of Grothendieck categories and left adjoints as morphisms. This choice of morphisms is broadly accepted as the correct notion of morphism between noncommutative schemes given by Grothendieck categories (see for example \cite{staffordvandenbergh01}, \cite{rosenberg98}).
	\item[$\mathsf{Grt}_\otimes$]  is the $2$-category of Grothendieck categories equipped with a tensorial structure and monoidal left adjoints as morphisms (see for example \cite{brandenburg2014tensor}, \cite{brandenburgchirvasitu14}). This choice seems very natural for a categorical approach to classical algebraic geometry, as given a scheme $X$, $\mathsf{Qcoh}(X)$ has a canonical monoidal structure. The study of this 2-category in the aforementioned \cite{brandenburg2014tensor} and \cite{brandenburgchirvasitu14} is motivated by the discussion of ``2-algebraic geometry'' from \cite{chirvasitujohnsonfreyd13} and by the reconstruction theorems for stacks.
	\item[$\grothflat$]  is the $2$-category of Grothendieck categories and left exact left adjoints as morphisms. From a naive categorical point of view, this choice is the one that takes more seriously the analogy with topoi (and of course geometric morphisms between them). These morphisms have been used in \cite{rosenberg98} and \cite{kontsevichrosenberg04} as the correct notion of a flat morphism between noncommutative schemes, motivating our choice of notation $\grothflat$. We will provide later some examples to illustrate this intuition.
	\end{itemize}
The $2$-category $\grothflat$ is the main object of study of this paper. More concretely, we show that $\grothflat$ can be endowed with a monoidal structure and we characterize the exponentiable objects therein. From an algebro-geometric point of view, this can be seen as a contribution to the understanding of exponentiable schemes or Hom-schemes when we restrict ourselves to the flat case.
 
\subsection{Monoidal structures on $2$-categories of noncommutative schemes} The 2-category $\Pres$ of locally presentable categories and cocontinuous functors can be endowed with a closed symmetric monoidal structure \cite[Ch. 5]{bird84} 
$$\boxtimes: \Pres \times \Pres \ra \Pres,$$
which is built upon Kelly's tensor product of $\alpha$-cocomplete small categories \cite{kelly82structures}, where $\alpha$ is any regular cardinal. As Kelly's tensor product is defined in the enriched setup, one can analogously recover a monoidal structure in the 2-category $\mathbb{Z}\mhyphen \Pres$ of locally presentable linear categories 
$$\boxtimes: \mathbb{Z}\mhyphen \Pres \times \mathbb{Z}\mhyphen \Pres \ra \mathbb{Z}\mhyphen \Pres$$
(see for example \cite[Cor 2.2.5]{chirvasitujohnsonfreyd13}). As Grothendieck categories are in particular locally presentable (see for example \cite[Prop 3.4.16]{borceuxhandbook3}), it is natural to wonder whether, given a 2-category of Grothendieck categories, the tensor product of locally presentable categories endows it with a monoidal structure.

\subsubsection{Monoidal structure on $\groth$}
As mentioned above, Grothendieck categories are precisely the linear Grothendieck topoi. Using this perspective, a tensor product of Grothendieck categories $\boxtimes$ was introduced in \cite{lowenramosgonzalezshoikhet18}. In particular, $\boxtimes$ can be seen as a linear parallel to the product of Grothendieck topoi \cite{pitts85} or to the product of $\infty$-topoi \cite{lurie09, lurie07}. This tensor product induces a monoidal structure
$$\boxtimes: \groth \times \groth \ra \groth$$ 
in $\groth$, which, contrary to the product of topoi or $\infty$-topoi, is a honest monoidal structure and not a cartesian product.
In addition, $\boxtimes$ can be seen to coincide with the tensor product of locally presentable linear categories, which explains our choice to express both with the same symbol.
 
\subsubsection{Monoidal structure on $\groth_\otimes$}
In \cite[\S 5]{brandenburg2020bicategorical} it is shown that the tensor product of locally finitely presented linear categories 
$$\boxtimes: (\mathbb{Z}\mhyphen\mathsf{Pres})_{\omega} \times (\mathbb{Z}\mhyphen\mathsf{Pres})_{\omega} \ra (\mathbb{Z}\mhyphen\mathsf{Pres})_{\omega}$$ respects tensorial structures, inducing a monoidal structure in the category of locally finitely presented tensor linear categories 
$$\boxtimes: (\mathbb{Z}\mhyphen\mathsf{Pres})_{\omega,\otimes} \times (\mathbb{Z}\mhyphen\mathsf{Pres})_{\omega,\otimes} \ra (\mathbb{Z}\mhyphen\mathsf{Pres})_{\omega,\otimes}.$$
It is most likely that an analogous argument will work for higher cardinalities, allowing us to obtain a monoidal structure 
$$\boxtimes: (\mathbb{Z}\mhyphen\mathsf{Pres})_{\otimes} \times (\mathbb{Z}\mhyphen\mathsf{Pres})_{\otimes} \ra (\mathbb{Z}\mhyphen\mathsf{Pres})_{\otimes}$$
for locally presentable tensor linear categories in general. As we know that $\boxtimes$ restricts to Grothendieck categories, in particular this would allow to obtain a monoidal structure on $\groth_\otimes$ as well.

From \cite[Thm A]{brandenburg20}, one can deduce that not only the functor
$$\Qcoh: \Sch^\op_{\qcqs} \ra \groth_\otimes$$ 
assigning to each quasi-compact quasi-separated scheme $X$ its Grothendieck tensor category of quasi-coherent sheaves, but also its composition with the forgetful functor
$$\Qcoh: \Sch^\op_{\qcqs} \ra \groth$$
are monoidal with respect to the product of schemes in $\Sch_{\qcqs}$ and the tensor product $\boxtimes$ of Grothendieck categories in $\grothflat$ and $\groth$. Consequently, in the setup of noncommutative algebraic geometry, $\boxtimes$ can be seen as \textit{the right notion} of tensor product of noncommutative schemes.  

\subsubsection{Monoidal structure on $\grothflat$}
As explained above, the tensor product $\boxtimes$ of locally presentable linear categories restricts to Grothendieck categories. Therefore, given $\cc, \cd \in \grothflat$, we have an assignment
$$(\cc,\cd) \mapsto \cc \boxtimes \cd \in \grothflat.$$
However, to show that this assignment is actually functorial in $\grothflat$ is far from trivial. The first half of our paper will be devoted to show that this is indeed the case, providing a monoidal structure 
$$\boxtimes: \grothflat \times \grothflat \ra \grothflat$$
in $\grothflat$. This can be seen as a noncommutative parallel of the fact that the product of flat morphisms of schemes is again flat \cite[Cor 2.1.7]{egaiv2}.

\subsection{Exponentiability} 
When a ($2$-)category $\mathsf{G}$ of Grothendieck categories (think of any of the examples introduced above) is equipped with a monoidal structure $\boxtimes$, it is a natural question to wonder whether an object $\cd$ is exponentiable with respect to $\boxtimes$, that is, whether for every object $\cb$ there exist a notion of internal hom $\cb^\cd$ such that \[\mathsf{G}(\ca \boxtimes \cc, \cb) \simeq \mathsf{G}(\ca, \cb^\cc), \]
naturally in $\ca$. As many readers may know, this is precisely the same as requesting a right adjoint for the functor $- \boxtimes \cc$. This question is not only natural, but also useful, as it provides an internal object of morphisms, from which one can then recover the external hom via the formula \[ \mathsf{G}(\cc, \cb) \simeq \mathsf{G}(1, \cb^\cc),\] where $1$ is the unit of the monoidal structure. An object of morphisms is a much better device than an external hom, as it is a Grothendieck category (in our case) and we have a whole cohomological machinery designed for those. From the perspective of noncommutative algebraic geometry, the exponentiable Grothendieck categories in $\mathsf{G}$ can serve as a formal replacement for the exponentiable schemes (in the suitable category of schemes of which $\mathsf{G}$ aims to provide a noncommutative version). In particular, the study of exponentiable objects in such a category $\mathsf{G}$ could shed some light in the grasp of exponentiable schemes, from which we do not have a full understanding.

\subsubsection{Exponentiable objects in $\groth$}
While it is well-known that categories of modules are exponentiable in $\groth$, there are also known counterexamples of exponentiability (see, for example, \cite[Rem 6.5]{positselskirosicky17}). The full characterization of exponentiable objects in $\groth$ is part of an ongoing joint research project of the second named author with Wendy Lowen and Michel Van den Bergh. 
\subsubsection{Exponentiable objects in $\groth_\otimes$}
To the best of our knowledge, the exponentiability of objects in $\groth_\otimes$ has not been analysed as such in the literature. Nonetheless, given the fact that the functor
$$\Qcoh: \Sch_{\qcqs}^\op \ra \groth_\otimes$$
is fully faithful (see \cite[Thm 3.4.3]{brandenburgchirvasitu14}), in the essential image of the functor, the question can be reduced to the characterization of exponentiability in $\Sch_{\qcqs}$, that is, to the study of the representability of the Hom-scheme functor. In this situation, it is known that, under suitable assumptions on the schemes involved (see \cite[\S 4c]{grothendieck95}), the internal scheme of morphisms $Y^X$ exists when $X$ is projective and flat and $Y$ is quasi-projective. Nevertheless, the Hom-scheme is not representable in $\Sch_{\qcqs}$ in general.
 
\subsubsection{Exponentiable objects in $\grothflat$} 
The second part of this paper will be devoted to establishing a characterization of the exponentiable objects in $\grothflat$. A characterization of the exponentiable objects in the category of Grothendieck topoi is provided in \cite{johnstonejoyal82} and, using similar techniques, a characterization of exponentiable $\infty$-topoi is given in \cite{anellejay18}. The fact that $\grothflat$ can be seen as ``the category of linear topoi'', allows us to follow a parallel strategy. 

\subsection{Structure of the paper} This paper provides a systematic study of the monoidal structure $\boxtimes$ on $\grothflat$ and describes its exponentiable objects. Besides the technical result per se, we see this as a contribution to the debate about adequate categories of (non)commutative schemes, stressing on the fact that the right choice of morphism is a key aspect of devising such a category. The paper is structured as follows:
\begin{itemize}
	\item[\S 2] is a soft introduction to our $2$-category of interest: $\grothflat$ (\S \ref{subsec:enrichedtopoi}). Besides recalling the main definitions and fixing the notations, the main purpose of this section is to convince the reader that $\grothflat$ can simulate flat algebraic geometry via a collection of examples (\Cref{contextflatalggeo}). Together with the first subsection of \S3, this is the only non-original part of the paper.
	\item[\S 3] We show that the \textbf{monoidal structure} $\boxtimes$ on $\mathsf{Grt}$ defined in \cite{lowenramosgonzalezshoikhet18} (\S\ref{subsecpresentations}) nicely restricts to $\grothflat$ (\S\ref{subsecmonoidal}). Notice that while this is a trivial task on the level of objects, it is a highly non-trivial task on the level of morphisms. We also introduce the problem of exponentiability (\S\ref{subsecexponentiable}) and we give a first and very partial result as a kind of motivational example: we show that categories of linear presheaves $\Mod(\AAA)$ are exponentiable (\Cref{propmodulesareexponentiable}).
	\item[\S 4] We study the properties of the forgetful functor $\mathbb{U}: \grothflat^{\op} \ra \Cat_k$ and we show that it is representable (\Cref{forgetfulrep}). The Grothendieck category that represents this forgetful functor will be denoted by $\Mod(k)[\obj]$.
	\item[\S 5] We introduce and study \textbf{(quasi-)injective Grothendieck categories} (\S\ref{qinjgcats}), \textbf{continuous linear categories} (\S\ref{contcat}) and then connect the two concepts (\S\ref{contvsqinj}). These will be relevant technical tools for our main theorem.
	\item[\S 6] Contains our \textbf{main theorem}:	
						\begin{thm*}[\Cref{maintheorem}]
						A Grothendieck category is exponentiable in $\grothflat$ if an only if it is continuous. In particular, every finitely presentable Grothendieck category is exponentiable.
						\end{thm*} 
	\item[\S 7] is a collection of examples and instances of our main theorem. The most relevant is the following proposition:
						\begin{thm*}[\Cref{qcohexponentiable}]
						Let $X$ be a quasi-compact quasi-separated scheme over $k$, then $\mathsf{Qcoh}(X)$ is exponentiable.
						\end{thm*}
	\noindent The rest of the section is dedicated to an indepth analysis of concrete examples where exponentiations can be computed more or less explicitely. 
	\end{itemize}

\section{Preliminaries}\label{secpreliminaries}
The main scope of this section is to introduce the $2$-category $\grothflat$ (\Cref{def:grothflat}), which is our main object of study throughout the paper, and relate it to a $2$-category $\mathsf{Site}_{\text{flat}}$ of linear sites (\Cref{def:siteflat}). The last subsection shows that $\grothflat$ is a good framework to accommodate (noncommutative) flat algebraic geometry.

\begin{notat} We fix $k$ a commutative ring. Throughout the rest of the paper we will work enriched over the category $\Mod(k)$ of (right) $k$-modules. The term \emph{$k$-linear} is frequently used in a categorical framework as a synonym of \emph{$\Mod(k)$-enriched}. We will follow this convention along the text.
\end{notat}

\subsection{The theory of $k$-linear sheaves}\label{subsec:enrichedtopoi}

The category $\Mod(k)$ of (right) $k$-modules is a locally presentable symmetric monoidal closed category which is also regular in the sense of Barr \cite{barr71}. These properties guarantee the existence of a well-behaved theory of $k$-linear sheaves \cite{borceuxquinteiro96}. More concretely, given a small $k$-linear category $\AAA$, there is a one to one correspondence between left exact reflections of $k$-linear presheaves on $\AAA$ and $k$-linear Grothendieck topologies on $\AAA$ \cite[Thm 1.5]{borceuxquinteiro96}. Therefore, we have a $k$-linear parallel to topos theory. In this section, we want to provide a short overview of the basic elements of this theory. 

\begin{notat}
Let $\AAA$ be a small $k$-linear category. The category of $k$-linear presheaves $k \mhyphen[\AAA^{\op},\Mod(k)]$ is usually denoted in the literature by $\Mod(\AAA)$  and referred to as the \emph{category of (right) $\AAA$-modules}. We will also adopt this notation.
\end{notat}

The Gabriel-Popescu theorem \cite{gabrielpopescu64} characterizes the localizations (i.e. the left exact reflections) of categories of modules. These are precisely the Grothendieck $k$-linear categories. We recall the definition.

\begin{defn}
	A \emph{Grothendieck $k$-linear category} is a $k$-linear abelian category with small colimits, exact filtered colimits and a generator. 
\end{defn}

On the other hand, following the theory of enriched sheaves from \cite{borceuxquinteiro96}, we can also see the Grothendieck $k$-linear categories as the categories of $k$-linear sheaves on $k$-linear sites. We recall now the basic definitions of the theory of $k$-linear Grothendieck topologies and sites. For an in-depth analysis, we point the reader to \cite[\S 2]{lowen16} and \cite[\S 2]{ramosgonzalez18}. 

Let $\AAA$ be a small $k$-linear category.
\begin{defn}
	Given an object $A \in \AAA$, a \emph{sieve on $A$} is a subobject $R$ of the representable module $\AAA(-,A)$ on $A$ in $\Mod(\AAA)$. Given a sieve $R$ on $A$ and a morphism $f:B \ra A$ in $\AAA$, the \emph{pullback sieve} $f^{-1}R \subseteq \AAA(-,B)$ is the sieve on $B$ given by $f^{-1}R(C) \coloneqq \{g: C \ra B \,\,|\,\, f\circ g \in R(C)\} \subseteq \AAA(C,B)$.
\end{defn}
 
\begin{defn}
	A \emph{$k$-linear Grothendieck topology} $\tau$ on $\AAA$ consists of, for each $A \in \AAA$, a family of sieves $\tau(A)$, called \emph{covering sieves}, on $A$ satisfying the following three axioms:
	\begin{itemize}
		\item For every $A \in \AAA$, $\AAA(-,A) \in \tau(A)$;
		\item For every morphism $f: B \ra A$ in $\AAA$, and every covering sieve $R \in \tau(A)$ the pullback sieve $f^{-1}R \subseteq \AAA(-,B)$ belongs to $\tau(B)$;
		\item Given a covering sieve $R \in \tau(A)$ and a sieve $S$ on $A$, if for every $B \in \AAA$ and every $f:B \ra A$ in $R(B)$ the sieve $f^{-1}S$ belongs to $\tau(B)$, then $S$ belongs to $\tau(A)$.
	\end{itemize}
 	A pair $(\AAA,\tau)$ of a small $k$-linear category $\AAA$ endowed with a $k$-linear Grothendieck topology $\tau$ is called a \emph{$k$-linear site}.
\end{defn}

\begin{defn}
	A family of morphisms $\cf = \{a_i: A_i \ra A\}_{i \in I}$ in a $k$-linear site $(\AAA,\tau)$ is \emph{covering} if the smallest $k$-linear sieve $S$ that contains $\cf$ belongs to $\tau(A)$. 
\end{defn}

\begin{defn}
	A \emph{$k$-linear sheaf} $F$ on $(\AAA, \tau)$ is a $k$-linear presheaf $F \in \Mod(\AAA)$ such that the restriction
	\begin{equation*}
	F(A) \simeq \Mod(\AAA)(\AAA(-,A), F) \ra \Mod(\AAA)(R, F)
	\end{equation*}
	is an isomorphism for all $A \in \AAA$ and all covering sieve $R \in \tau(A)$. 
	We denote by $$\Sh(\AAA, \tau) \subseteq \Mod(\AAA)$$ the full subcategory of $k$-linear sheaves, and by 
	$$\#: \Mod(\AAA) \ra \Sh(\AAA,\tau)$$
	the sheafification ($k$-linear) functor.
\end{defn}

In what follows, we provide the $k$-linear analogues of the (2-)category of topoi and the (2-)category of sites.

\subsubsection{The topos-like 2-category of $k$-linear Grothendieck categories} \label{toposlikegrtflat}

We proceed to describe the $k$-linear parallel of the $2$-category of Grothendieck topoi, that will denote by $k \mhyphen\grothflat$. 

The objects of $k \mhyphen\grothflat$ are the Grothendieck $k$-linear categories. The $k$-linear analogue of the geometric morphisms is given as follows.  
\begin{defn}\label{defgeometricmorphism}
	Consider two $k$-linear Grothendieck categories $\cc, \cd$. A \emph{$k$-linear flat morphism} $F: \cc \ra \cd$ is an adjunction $F^*: \cd \rightleftarrows \cc : F_*$ of $k$-linear functors such that the left adjoint $F^*$ is exact.
\end{defn}

\begin{rem}
	Observe that, by the Special Adjoint Functor Theorem, we have that giving a flat morphism $F:\cc \ra \cd$ is equivalent to giving a colimit preserving left exact functor $F^*:\cd \ra \cc$.
\end{rem}
 	
In addition, we have the $k$-linear version of morphisms between geometric morphisms.
\begin{defn}\label{defgeometricnaturaltransformation}
	Given $F,G: \cc \ra \cd$ two $k$-linear flat morphisms, a morphism $F \Rightarrow G$ is a $k$-linear natural transformation $F^* \Rightarrow G^*$ between the left adjoints. 
\end{defn}

\begin{defn}\label{def:grothflat}
	The 2-category of $k$-linear Grothendieck categories $k \mhyphen\grothflat$ is the $2$-category with objects the $k$-linear Grothendieck categories, $k$-linear flat morphisms as $1$-cells and $k$-linear natural transformations between the left adjoints as $2$-cells.
\end{defn}

\begin{rem}[Why ``flat morphisms'' and not ``geometric morphisms''?]
	The morphisms defined in \Cref{defgeometricmorphism} are precisely the flat morhpisms between Gro\-thendieck categories in the sense of \cite{rosenberg98} or \cite{kontsevichrosenberg04}. As explained in the introduction, these can be thought as the noncommutative analogue of the classical flat morphisms of schemes. We stick to this terminology to emphasize this fact.
\end{rem}

\subsubsection{The 2-category of $k$-linear sites} \label{2catsites}
We now proceed to describe the $2$-category of $k$-linear sites, that we will denote by $k \mhyphen\siteflat$. 

The objects of this 2-category are the $k$-linear sites. 

The $k$-linear analogue of a morphism of sites is given as follows.
\begin{defn}\label{def:morphismsites}
	Consider $k$-linear sites $(\AAA,\tau_{\AAA})$ and $(\BBB,\tau_{\BBB})$. A $k$-linear functor $f: \AAA \ra \BBB$ is a \emph{morphism of $k$-linear sites}  if:
	\begin{itemize} 
		\item It is \emph{covering flat}, i.e. the functor 
		$$\Mod(\AAA) \ra \Sh(\BBB) : M \mapsto (\#_{\BBB} \Lan_{Y_{\AAA}}(Y_\BBB f))(M)$$ preserves finite limits and
		\item it is \emph{cover-preserving}  (sometimes also called \emph{continuous}), i.e. for every covering family $\{a_i: A_i \ra A\}_{i \in I}$ in $\tau_{\AAA}$, the family $\{f(a_i):f(A_i) \ra f(A)\}_{i\in I}$ is covering in $\tau_{\BBB}$.
	\end{itemize}
	If this is the case, it can be easily shown that $f$ induces a flat morphism between the categories of sheaves. More concretely, if a morphism $f:\AAA \ra \BBB$ is cover-preserving, $f^*: \Mod(\BBB) \ra \Mod(\AAA)$ restricts to a functor $f_s: \Sh(\BBB,\tau_{\BBB}) \ra \Sh(\AAA,\tau_{\AAA})$, which has a left adjoint $f^s: \Sh(\AAA,\tau_{\AAA}) \ra \Sh(\BBB,\tau_{\BBB})$. If in addition $f$ is covering flat, then the adjunction $f^s: \Sh(\AAA,\tau_{\AAA}) \rightleftarrows \Sh(\BBB,\tau_{\BBB}): f_s$ is a flat morphism. 
\end{defn}

\begin{rem} The notation $(f^s,f_s)$ for the adjunction between the sheaf categories induced by $f$ is borrowed from \cite{SGA4-1}. Though this notation is not standard in the topos theory community, we believe that in the setup of \Cref{def:morphismsites} is the most suitable choice in order to make a clear disctintion with the adjunction $(f^*,f_*)$ between the presheaf categories induced by $f$.
\end{rem}

The $2$-cells will be given by the natural choice:
\begin{defn}
	Given $f,g: (\AAA,\tau_{\AAA}) \ra (\BBB,\tau_{\BBB})$ two morphisms of sites, a morphism $f \Rightarrow g$ is just a $k$-linear natural transformation $f \Rightarrow g$. 
\end{defn}

\begin{defn}\label{def:siteflat}
The 2-category of $k$-linear sites $k \mhyphen\siteflat$ is the $2$-category with objects the $k$-linear sites, $k$-linear morphisms of sites as $1$-cells and $k$-linear natural transformations as $2$-cells.
\end{defn}

\begin{rem}[Why the subindex ``flat''?]
	As in the case of Grothendieck categories, depending on the context one may want to consider different $2$-categories of $k$-linear sites. In our setup, parallely to classical topos theory, morphisms of sites allow to define a pseudofunctor $k \mhyphen\siteflat \ra k \mhyphen\grothflat^{\op}$ (see \cite[\S4]{ramosgonzalez18}). However, if we were to be more interested in working with the $2$-category $\groth$, it is then more natural to consider the 2-category of $k$-linear sites with 1-cells just the cover-preserving morphisms. Indeed, if we call this 2-category $k \mhyphen\site$, we then obtain a pseudofunctor $k \mhyphen\site \ra k \mhyphen\groth$ (see \cite[\S4]{ramosgonzalez18}).
\end{rem}

The following result relating flat morphisms and morphisms is the $k$-linear parallel of \cite[Cor C.2.3.9]{elephant2}.

\begin{prop}\label{propgeometricmorphismssites}
	Given a Grothendieck category $\ca$ and a $k$-linear site $(\BBB,\tau_{\BBB})$, we have that
	\begin{equation*}
	k \mhyphen\grothflat(\ca,\Sh(\BBB,\tau_{\BBB})) = k \mhyphen\siteflat((\BBB,\tau_{\BBB}),\ca),
	\end{equation*}
	where $\ca$ is considered as a site endowed with its canonical topology.

\end{prop}

 \subsection{Examples from flat algebraic geometry} \label{contextflatalggeo}
In the framework of noncommutative algebraic geometry, flat morphisms between abelian categories are used in \cite{rosenberg98} and \cite{kontsevichrosenberg04} as the correct noncommutative version of the classical flat morphisms between schemes. In this subsection we want to explain the intuition behind this idea through a series of examples, illustrating the relevance of $k \mhyphen\grothflat$ in the study of flat algebraic geometry.  The reader that is somewhat new to flatness and is looking for a geometric interpretation of such concept might find \cite{6789} interesting.

\begin{exa}
Consider the ringed space $(\Spec(k),\mathcal{O}_{\Spec(k)})$. Given a $k$-ringed space $(X, \mathcal{O}_X)$ (i.e. a morphism $(X, \mathcal{O}_X) \ra (\Spec(k),\mathcal{O}_{\Spec(k)})$ of ringed spaces), its category of sheaves of $\mathcal{O}_X$-modules is a Grothendieck $k$-linear category \cite[Tag 01AH]{stacksproject}. If we denote by $k \mhyphen \Ringspaces$ the category of $k$-ringed spaces, we have a functor 
\begin{equation}
\mathsf{Mod}: k\mhyphen\Ringspaces^{\op} \to k \mhyphen \groth
\end{equation}
that sends each $k$-ringed space $(X, \mathcal{O}_X)$ to its category of presheaves of $\mathcal{O}_X$-modules $\Mod(X)$ and to every morphism $f:(X,\mathcal{O}_X) \ra (Y,\mathcal{O}_Y)$ of $k$-ringed spaces to the pullback functor $f^*:\Mod(Y) \ra \Mod(X)$, which is a $k$-linear functor. 
If we now consider the non-full subcategory $k\mhyphen \Ringspaces_{\fflat} \subseteq k \mhyphen \Ringspaces$ of $k$-ringed spaces with flat morphisms, we obtain a functor 
\begin{equation}
\mathsf{Mod}: k\mhyphen\Ringspaces_{\fflat} \to k \mhyphen\grothflat
\end{equation}
that sends each $k$-ringed space $(X, \mathcal{O}_X)$ to its category of presheaves of $\mathcal{O}_X$-modules $\Mod(X)$ and to every morphism $f:(X,\mathcal{O}_X) \ra (Y,\mathcal{O}_Y)$ to the pullback-pushforward $k$-linear adjunction $f^*:\Mod(Y) \rightleftarrows \Mod(X): f_*$, which is a $k$-linear flat morphism as a direct consequence of the flatness of $f$ (see \cite[Tag 02N4]{stacksproject}).
\end{exa}

\begin{exa}
If we restrict the setup of the previous remark by taking $k$-schemes (i.e. schemes relative over $(\Spec(k),\mathcal{O}_{\Spec(k)})$) instead of $k$-linear ringed spaces, and quasi-coherent sheaves of modules instead of all presheaves of modules, we obtain the following.
Given a $k$-scheme $X$, the category of quasi-coherent sheaves $\text{Qcoh}(X)$ is a $k$-linear Grothendieck category \cite[Tag 077P]{stacksproject}. Denote by $k \mhyphen \Sch_{\qcqs}$ the category of $k$-schemes with quasi-compact quasi-separated morphisms of $k$-schemes. We have a functor
\begin{equation}
	\Qcoh: k\mhyphen \Sch_{\qcqs}^{\op} \ra k \mhyphen \groth
\end{equation}
that sends each $k$-scheme $X$ to its category of quasi-coherent sheaves $\Qcoh(X)$ and each quasi-compact quasi-separated morphism of $k$-schemes $f:X \ra Y$ to the pullback functor $f^*: \Qcoh(Y) \rightleftarrows \Qcoh(X): f_*$ (see \cite[Tag 01LC]{stacksproject}), which is a $k$-linear functor.

If we now consider the non-full subcategory $k \mhyphen \Sch_{\qcqs,\fflat} \subseteq k \mhyphen\Sch_{\qcqs}$ of $k$-schemes with quasi-compact quasi-separated flat morphisms, we have a functor
\begin{equation}
	\Qcoh: k \mhyphen \Sch_{\qcqs,\fflat} \ra k \mhyphen \grothflat
\end{equation}
that sends each $k$-scheme $X$ to its category of quasi-coherent sheaves $\Qcoh(X)$ and each quasi-compact quasi-separated flat morphism of schemes $f:X \ra Y$ to the pullback-pushforward $k$-linear adjunction $f^*: \Qcoh(Y) \rightleftarrows \Qcoh(X): f_*$, which is a $k$-linear flat morphism as a direct consequence of the flatness of $f$.
\end{exa}

\section{The tensor product of Grothendieck categories}

The tensor product $\boxtimes_k$ of Grothendieck $k$-linear categories introduced in \cite{lowenramosgonzalezshoikhet18} is shown in \cite{ramosgonzalez20} to endow $k \mhyphen \groth$ with a symmetric monoidal structure. In this section, after revising in \S\ref{subsecpresentations} the different presentations of $\boxtimes_k$, we show in \S\ref{subsecmonoidal} that $\boxtimes_k$ also endowes $k\mhyphen \grothflat$, and not only $k \mhyphen \groth$, with a symmetric monoidal structure. To conclude, we introduce in \S\ref{subsecexponentiable} the main goal of the paper, that is, the characterization of the exponentiable objects in $k \mhyphen \grothflat$ with respect to the monoidal structure given by $\boxtimes_k$. 

\begin{notat}
In order to ligthen the notations, we will omit the prefix ``$k \mhyphen$'' from now on. Unless explicitly stated, all categories and functors considered will be $k$-linear. In particular, we write $\Pres = k$-$\Pres$, $\groth= k$-$\groth$ and $\grothflat = k$-$\grothflat$.
\end{notat}

\subsection{The tensor product of Grothendieck categories: Presentations}\label{subsecpresentations}

As explained in \S\ref{secpreliminaries}, the theory of linear sheaves provides us with a one to one correspondence between linear Grothendieck topologies on a small linear category $\AAA$ and left exact reflections of $\Mod(\AAA)$, allowing us to think of Grothendieck categories as a linear parallel to Grothendieck topoi. In \cite{lowenramosgonzalezshoikhet18} a tensor product of Grothendieck categories $\boxtimes$ is defined in terms of both a tensor product of linear sites and a tensor product of left exact reflections of module categories, which are shown to agree via the aforementioned one to one correspondence. On the other hand, Grothendieck categories are in particular instances of locally presentable linear categories and $\boxtimes$ is shown in loc. cit. to coincide with the tensor product of these latter, providing us with a formal presentation of $\boxtimes$, independent of a choice of representative. 

It is important to remark that, though the tensor product of Grothendieck categories is a linear version of the product of Grothen\-dieck topoi, it differs from this latter on the fact that it is not a categorical product in $\grothflat$, but a honest monoidal structure, as it will become evident from its presentations.

\subsubsection{The tensor product in terms of presentations}
In this subsection we describe the tensor product of linear sites and the tensor product of left exact reflections of module categories from \cite{lowenramosgonzalezshoikhet18}. These tensor products give rise to a well-defined tensor product of Grothendieck categories. 

\begin{defn}
	Let $(\AAA,\tau_{\AAA}), (\BBB,\tau_{\BBB})$ be linear sites. Given $A \in \AAA$, $B \in \BBB$ and covering sieves $R \in \tau_{\AAA}(A)$ and $S \in \tau_{\BBB}(B)$, consider the canonical morphism
	$$\phi_{R,S}: R \otimes S \ra  \AAA(-,A) \otimes \BBB(-,B) = \AAA \otimes \BBB(-, (A,B)).$$
	The \emph{tensor product sieve} of $R$ and $S$ is defined as $R \boxtimes S = \Img(\phi_{R, S})$. We define the \emph{tensor product topology} $\tau_{\AAA} \boxtimes \tau_{\BBB}$ on $\AAA \otimes \BBB$ to be the smallest topology containing the cover system $\crr = \{ R \boxtimes S \,\, |\,\, R \in \tau_{\AAA}, S \in \tau_{\BBB} \}$. We define the \emph{tensor product of linear sites} as $(\AAA,\tau_{\AAA}) \boxtimes (\BBB,\tau_{\BBB}) \coloneqq (\AAA \otimes \BBB, \tau_{\AAA} \boxtimes \tau_{\BBB})$.
\end{defn}

\begin{defn}\label{defntpgrothcatlocalisations}
	Given two left exact reflections of module categories $\ca \subseteq \Mod(\AAA)$, $\cb \subseteq \Mod(\BBB)$ respectively, we define the \emph{tensor product of left exact reflections} $(\ca \subseteq \Mod(\AAA))\boxtimes (\cb \subseteq \Mod(\BBB))$ as the full subcategory of $\Mod(\AAA \otimes \BBB)$ spanned by the objects
	\begin{equation*}
	\{ F \in \Mod(\AAA\otimes\BBB) \,\,|\,\, F(A,-) \in \cb \,\,\forall A \in \AAA, F(-,B) \in \ca \,\,\forall B \in \BBB\}.
	\end{equation*}
	This is indeed a left exact reflection of $\Mod(\AAA \otimes \BBB)$. 
\end{defn}

\begin{prop}[{\cite[Cor 2.16, Cor 2.18, Prop 4.1]{lowenramosgonzalezshoikhet18}}] \label{proptensorproductlocalisations}
Consider two Grothendieck categories $\ca, \cb$ and choose linear sites $(\AAA,\tau_{\AAA})$ and $(\BBB,\tau_{\BBB})$ such that $\ca \simeq \Sh(\AAA,\tau_{\AAA}) \subseteq \Mod(\AAA), \cb \simeq \Sh(\BBB, \tau_{\BBB}) \subseteq \Mod(\BBB)$. Then we have that
\begin{equation*}
	(\ca \subseteq \Mod(\AAA))\boxtimes (\cb \subseteq \Mod(\BBB)) \simeq \Sh((\AAA,\tau_{\AAA})\boxtimes(\BBB, \tau_\BBB)).
\end{equation*} 
In particular, this does not depend on the representations of $\ca$ and $\cb$ chosen, that is, given other linear sites $(\AAA',\tau_{\AAA'})$ and $(\BBB',\tau_{\BBB'})$ such that $\ca \simeq \Sh(\AAA',\tau_{\AAA'}) \subseteq \Mod(\AAA'), \cb \simeq \Sh(\BBB', \tau_{\BBB'}) \subseteq \Mod(\BBB')$, one has that 
\begin{equation*}
	\begin{aligned}
		(\ca \subseteq \Mod(\AAA))\boxtimes (\cb \subseteq \Mod(\BBB))&\simeq (\ca \subseteq \Mod(\AAA'))\boxtimes (\cb \subseteq \Mod(\BBB'))\\ 
		&\simeq \Sh((\AAA',\tau_{\AAA'})\boxtimes(\BBB', \tau_{\BBB'})) \simeq \Sh((\AAA,\tau_{\AAA})\boxtimes(\BBB, \tau_\BBB))
	\end{aligned} 
\end{equation*}
\end{prop}
Therefore, we can define a \emph{tensor product of Grothendieck categories} $\ca \boxtimes \cb$ up to equivalence of categories, given by taking presentations of $\ca$ and $\cb$ either as categories of linear sheaves or as left exact reflections and computing the corresponding tensor products of presentations \cite[Def 4.2]{lowenramosgonzalezshoikhet18}.  

\subsubsection{Tensor product à la Bird-Kelly: a formal description of the tensor product}
It is well-known that Grothendieck categories are in particular instances of locally presentable linear categories \cite[Prop 3.4.16]{borceuxhandbook3}. More concretely, we have that $\groth$ is a full sub-$2$-category of $\Pres$, where $\Pres$ is the $2$-category of locally presentable linear categories with $1$-cells the colimit preserving linear functors $\Cocont(-,-)$ and $2$-cells the linear natural transformations. It was shown by Bird in \cite{bird84} that (the non-linear version of) the $2$-category $\Pres$ has a monoidal structure given by \[\ca \boxtimes \cb = \Cont(\ca^\circ, \cb),\] which, together with the internal hom, equips (the non-linear version of) $\Pres$ with a monoidal closed structure. The linear case can be found, for example, in \cite{brandenburgchirvasitujohnsonfreyd15}. 

In parallel with the situation in the classical topos theory \cite[Thm 2.3]{pitts85}, the tensor product of Grothendieck categories $\boxtimes$ as defined in the previous section coincides with the tensor product of locally presentable linear categories \cite[Thm 5.4]{lowenramosgonzalezshoikhet18}. Consequently, $\groth$ is a full monoidal sub-2-category of $\Pres$, and in particular, the tensor product of Grothendieck categories $\boxtimes$ endowes $\groth$ with a monoidal structure.

This observation immediately allows us to obtain many useful properties of $\boxtimes$. We list some of them:

\begin{itemize}
	\item The tensor product of Grothendieck categories inherits the universal property of the tensor product of locally presentable categories, namely, given Grothendieck categories $\ca,\cb$ and $\cc$, we have that:
	\begin{equation}
	\Cocont(\ca \boxtimes \cb,\cc) = \Cocont(\ca,\Cocont(\cb,\cc)).
	\end{equation}
	\item The tensor product of Grothendieck categories is symmetric.
	\item As the unit of the tensor product of locally presentable $k$-linear categories is given by $\Mod(k)$ \cite{brandenburgchirvasitujohnsonfreyd15}, and $\Mod(k)$ is a Grothendieck $k$-linear category, we obtain that the unit of the tensor product of Grothendieck categories is also $\Mod(k)$.
\end{itemize}

\subsection{The monoidal structure on $\grothflat$} \label{subsecmonoidal}
The $2$-category $\groth$ is a monoidal full sub-2-category of $\Pres$, namely, the functor $-\boxtimes -: \Pres \times \Pres \ra \Pres$ restricts to a functor $-\boxtimes - :\groth \times \groth \ra \groth$. However, it is not immediately obvious whether this functor further restricts to the non-full sub-2-category $\grothflat^\op$, as it is not clear a priori whether the tensor product (of left adjoints) of flat morphisms is again (a left adjoint of) a flat morphism. This subsection will be devoted to show that this is indeed the case, providing us with a symmetric monoidal structure in $\grothflat^\op$ and thus in $\grothflat$. 
In other words, we wish to show that the tensor product of Grothendieck categories $\boxtimes$ induces a monoidal structure 
$$\boxtimes:\grothflat \times \grothflat \ra \grothflat$$
such that, given flat morphisms $f:\ca \ra \cb$ and $g:\cc \ra \cd$, the left adjoint $(f \boxtimes g)^*: \ca \boxtimes \cc \ra \cb \boxtimes \cd$ coincides with the tensor product $f^* \boxtimes g^*$ of $f^*$ and $g^*$ seen as morphisms in $\groth$. 

We organize this section as follows:
\begin{enumerate}
	\item Given flat morphisms $f$ and $g$ as above, we provide an implicit construction of the adjunction $(f \boxtimes g)^* \dashv (f \boxtimes g)_*$, where $(f\boxtimes g)^* \coloneqq f^* \boxtimes g^*$;
	\item We show that this adjoint pair is indeed a flat morphism and hence belongs to $\grothflat$. In order to accomplish this, we introduce a more concrete presentation of $(f\boxtimes g)^*$ which allows to easily check the left exactness.
\end{enumerate}

\subsubsection{An implicit description of the adjunction $(f \boxtimes g)^* \dashv (f \boxtimes g)_*$}\label{descrightadjoint}

Let $f: \ca \to \cb$ and $g: \cc \to \cd$ be flat morphism. In this subsection we furnish an implicit description of the adjunction
$$(f\boxtimes g)^*: \cb \boxtimes \cd \ra \ca \boxtimes \cc: (f\boxtimes g)_*,$$
where $(f \boxtimes g)^*$ is the cocontinuous functor provided by the tensor product $f^* \boxtimes g^*$ in $\groth$. 

First observe that, given the symmetry of the tensor product, we are reduced to describe the adjunction $(1_\cg \boxtimes g)^* \dashv (1_\cg \boxtimes g)_*$ for any Grothendieck category $\cg$ and any flat morphism $g: \cc \ra \cd$. One can then compute the adjunctions of the form $(g \boxtimes 1_\cg)^* \dashv (g\boxtimes 1_\cg)_*$ by means of the symmetry of $\boxtimes$, which can be implicitly described, in terms of the presentations of $\boxtimes$ as categories of continuous functors, as follows
$$\Cont(\ca^\op, \cb) \xrightarrow{\simeq} \Cont(\cb^\op,\ca): s \mapsto L(s)^\op,$$
where $L(s)$ denotes the left adjoint of $s$, which exists by the Adjoint Functor Theorem. 

We begin with the description of $(1_\cg \boxtimes f)_*$. Given the universal property of $\boxtimes$ in $\groth$, it is not hard to see that for any Grothendieck category $\cg$, the functor $(1_\cg \boxtimes g)_*:\cg \boxtimes \cc \ra \cg \boxtimes \cd$ can be explicitly described as
\begin{equation}\label{diagrightadjoint}
		(1_\cg \boxtimes g)_*: \Cont(\cg^\op,\cc) \ra \Cont(\cg^\op,\cd): s \mapsto g_* \circ s.
\end{equation}
Observe that $(1_\cg \boxtimes g)_*$ is indeed well-defined and continuous. 
\begin{rem}\label{remtheothervariable}
	 Let $f: \ca \to \cb$ be another flat morphism. By means of the symmetry, it is easy to check that the diagram
	\begin{equation*}
	\begin{tikzcd}
		{\Cont(\cg^\op,\ca)} \arrow[rrr, "f_* \circ -" description] \arrow[d, "\simeq" description] &  &  & {\Cont(\cg^\op,\cb)} \arrow[d, "\simeq" description] \\
		{\Cont(\ca^\op,\cg)} \arrow[rrr, "- \circ (f^*)^\op" description]                           &  &  & {\Cont(\cb^\op,\cg)}                                
	\end{tikzcd}
	\end{equation*}
	is commutative. In particular, we have that $(f \boxtimes 1_\cg)_*$ is given by
	\begin{equation*}
		(f\boxtimes 1_\cg)_*: \Cont(\ca^\op,\cg) \ra \Cont(\cb^\op,\cg): s \mapsto s \circ (f^*)^\op
	\end{equation*}
	and therefore, $(f \boxtimes g)_*$ is given by
	\begin{equation*}
		(f\boxtimes g)_*: \Cont(\ca^\op,\cc) \ra \Cont(\cb^\op,\cd): s \mapsto g_* \circ s \circ (f^*)^\op
	\end{equation*}
	or expressed diagramatically, $s \in \Cont(\ca^\op,\cc)$ is sent to the composite
	\begin{equation}
		\begin{tikzcd}
			\ca^\circ \arrow[rrrr, "s" description]  &  &  &  & \cc \arrow[dd, "g_*" description] \\
			&  &  &  &                                 \\
			\cb^\circ \arrow[rrrr, "g_* \circ s \circ (f^*)^\op" description, dashed] \arrow[uu, "(f^*)^\op" description] &  &  &  & \cd,                            
		\end{tikzcd}
	\end{equation}
	which lies in $\Cont(\cb^\op,\cd)$.
\end{rem}

On the other hand, a description of the functor $(1_\cg \boxtimes g)^*: \cg \boxtimes \cd \ra \cg \boxtimes \cc$ making use of the presentations of the tensor product as continuous functors is not readily available. Inspired by \eqref{diagrightadjoint}, one is inclined to think that the assignation
\begin{equation*}
		s \mapsto g^* \circ s \;\;\text{for }s \in \Cont(\cg^\op,\cd)
\end{equation*}
is the natural candidate for $(1_\cg \boxtimes g)$. Unfortunately, one immediately observes that $g^* \circ s$ is not continuous in general. Consequently, this assignation only yields a functor
\begin{equation*}
		(1_\cg \cap g)^*: \Cont(\cg^\op,\cd) \ra \Cat(\cg^\op, \cc):s \mapsto g^* \circ s
\end{equation*}
In order to find the correct description of $(1_\cg \boxtimes g)^*$ we observe that $\Cont(\cg^{\op}, \cc)$ is reflective\footnote{The careful reader will have noticed that we switched from the notation $\text{Cat}(\cg^{\op}, \cc)$  to $\underline{\text{Cat}(\cg^{\op}, \cc)}$. This is to handle size issues correctly. $\underline{\text{Cat}(\cg^{\op}, \cc)}$ is the full subcategory of $\text{Cat}(\cg^{\op}, \cc)$ spanned by (co)small functors (in the sense of \cite[\S2]{DAY2007651}). Indeed, every continuous functor is (co)small, because $\ca^\circ$ has a codense cogenerator, so the inclusion $\Cont(\cg^{\op}, \cc) \hookrightarrow \text{Cat}(\cg^{\op}, \cc)$ factors throught $\underline{\text{Cat}(\cg^{\op}, \cc)}$. The advantage of $\underline{\text{Cat}(\cg^{\op}, \cc)}$ is that it is a locally small category, differently from its bigger brother, which can be locally large somethimes. Now, since $\Cont(\cg^{\op}, \cc)$ is still a Grothendieck category, it is in particular cototal, and by \cite[Thm 1]{WOOD1982538} a continuous functor out of a cototal category into a locally small one has a left adjoint.} in $\underline{\text{Cat}(\cg^{\op}, \cc)}$. 
Let us call $L_{\cg,\cc}: \underline{\text{Cat}(\cg^{\op}, \cc)} \to \Cont(\cg^{\op}, \cc)$ the reflector. We claim that $(1_\cg \boxtimes g)^*$ is given by the following composition\footnote{In this diagram we are implicitely assuming that $(1_\cg \cap g)^*(s)$ is a (co)small functor for every $s$. Indeed this is true, because $g^*,s$ are (co)accessible functors, and thus their composition is also (co)small.}
\begin{equation}\label{eqleftadjoint}
\begin{tikzcd}
	{\Cont(\cg^{\op}, \cd)} \arrow[rrdd, "(1_\cg \cap g)^*" description] \arrow[rrrr, "(1_\cg \boxtimes g)^*", dashed] &  &                                                            &  & {\Cont(\cg^{\op}, \cc)} \\
	&  &                                                            &  &                         \\
	&  & {\underline{\text{Cat}(\cg^{\op}, \cc)}} \arrow[rruu, "L_{\cg,\cc}" description] &  &                        
\end{tikzcd}	
\end{equation}
Indeed, given the adjoint pair $g^* \circ -: \underline{\Cat(\cg^\op,\cd)} \rightleftarrows \underline{\Cat(\cg^\op,\cc)}: g_* \circ -$, we obtain the following chain of equivalences
\begin{equation*}
\begin{aligned}
		\Cont(\cg^\op,\cc)(L_{\cg,\cc}(g^* \circ s), t) &\simeq \underline{\Cat(\cg^\op,\cc)}(g^* \circ s, t) \\
		&\simeq \underline{\Cat(\cg^\op,\cd)}(s, g_* \circ t)\\ 
		&\simeq \Cont(\cg^\op,\cd)(s, (1_\cg \boxtimes g)_*(t)),
\end{aligned}
\end{equation*}
which are natural in $s \in \Cont(\cg^\op,\cd)$ and $t \in \Cont(\cg^\op,\cc)$, proving that the diagram \eqref{eqleftadjoint} provides an implicit description of $(1_\cg \boxtimes g)^*$, as desired. 

In a nutshell, $(1_\cg \boxtimes g)^*(s)$ is the \textit{best approximation} of $g^*\circ s$ among continuous functors.
 
\begin{rem}
	From \Cref{remtheothervariable}, one may think that a natural candidate for the left adjoint $(f \boxtimes 1_\cg)^*$ for a flat morphism $f:\ca \ra \cb$ would be given by the formula
	\begin{equation}\label{eqdescriptionleftadjointfail}
		\Cont(\cb^\op,\cg) \ra \Cont(\ca^\op,\cg): L_{\ca,\cg}(s \mapsto s \circ (f_*)^\op).
	\end{equation}
	This would be indeed the case if the diagram 
	\begin{equation*}
		\begin{tikzcd}
			{\Cont(\cg^\op,\cb)} \arrow[d, "\simeq" description] \arrow[rr, "f^* \circ-" description] &  & {\underline{\Cat(\cg^\op,\ca)}} \arrow[rr, "{L_{\cg,\ca}}" description] &  & {\Cont(\cg^\op,\ca)} \arrow[d, "\simeq" description] \\
			{\Cont(\cb^\op,\cg)} \arrow[rr, "-\circ (f_*)^\op" description]                           &  & {\underline{\Cat(\ca^\op,\cg)}} \arrow[rr, "{L_{\ca,\cg}}" description] &  & {\Cont(\ca^\op,\cg)}                                
		\end{tikzcd}
	\end{equation*}
	was commutative. However, as an explicit description of the reflectors $L_{\ca,\cg}, L_{\cg,\ca}$ is not in general at our disposal, it is not obvious how to check the commutativity of this diagram. For this reason, we unfortunately cannot conclude from this particular argument that \eqref{eqdescriptionleftadjointfail} indeed describes $(f \boxtimes 1_\cg)^*$.
\end{rem}

\begin{rem}
	Lamentably, from the description of $(1_\cg \boxtimes g)^*$ we just provided, one cannot deduce immediately that it is a left exact functor, as we do not know a priori if that is the case for the reflector $L_{\cg,\cc}$. In order to prove the left exactness of $(1_\cg \boxtimes g)^*$, we will devote next subsection to provide a more concrete description of this functor, for which left exactness will follow easily. 
\end{rem}

\subsubsection{Proving that $(f\boxtimes g)^*$ is left exact: a more concrete description}
As mentioned above, the symmetry of $\boxtimes$ allows us to work in one variable. Let $g: \cc \to \cd$ be a flat morphism and $\cg$ a Grothendieck category. In this subsection we provide a new description of $(1_\cg \boxtimes g)^*$ that will allow us to prove its left exactness. The structure of the argument will go through the following three steps, that we will condense in three propositions:
\begin{enumerate}
	\item We describe $(1_\cg \boxtimes g)^*$ when $\cg$ is locally finitely presentable and prove it is left exact (\Cref{leftexactfinitelypresentable});
	\item We describe $(f \boxtimes 1_\cg)^*$ when $f$ is a certain flat embedding and prove it is left exact (\Cref{leftexactalphaembedding});
	\item We combine the two previous results to show that $(1_\cg \boxtimes g)^*$ is left exact for any flat morphism $g: \cc \ra \cd$ and any Grothendieck category $\cg$ (\Cref{leftexactgeneralcase}).
	\end{enumerate}

\begin{prop}\label{leftexactfinitelypresentable}
Let $\cg$ be a locally finitely presented Grothendieck category and let $g: \cc \to \cd$ be a flat morphism between Grothendieck categories. Then, the left adjoint $(1_\cg \boxtimes g )^*$ of the functor $(1_\cg \boxtimes g)_*$ is given by 
\begin{equation}\label{comm1}
\begin{tikzcd}
\Cont(\cg^\op,\cc) \arrow[d, "\simeq" description]      &  & \Cont(\cg^\op,\cd) \arrow[ll, dashed, "(1_\cg \boxtimes g)^*" description] \arrow[d, "\simeq" description ] \\
{\Lex(\cg_\omega^{\op},\cc)} &  & {\Lex(\cg_{\omega}^{\op},\cd).} \arrow[ll, "g^* \circ -" description]                                           
\end{tikzcd}
\end{equation}
In particular, $(1_\cg \boxtimes g)^*$ is left exact.
\begin{proof}
	By hypothesis $\cg$ is locally finitely presentable and thus $\cg \simeq \Lex(\cg_{\omega}^\op,\Mod(k))$. Moreover, we know from \S\ref{descrightadjoint} that $(1_\cg \boxtimes g)_*: \Cont(\cg^\op,\cc) \ra \Cont(\cg^\op,\cd)$ is given by $g_*\circ -$. One can then easily check that the following diagram is commutative
	\begin{equation}\label{comm2}
		\begin{tikzcd}
		{\Cont(\cg^{\op},\cc)} \arrow[d, "\simeq" description ] \arrow[rr, "(1_\cg \boxtimes g)_*" description] &  & {\Cont(\cg^{\op},\cd)} \arrow[d,"\simeq" description] \\
		{\Lex(\cg_{\omega}^{\op},\cc)} \arrow[rr, "g_* \circ -" description]                       &  & {\Lex(\cg_{\omega}^{\op},\cd).}                       
		\end{tikzcd}
	\end{equation}
	Using the adjunction $g^*: \cd \rightleftarrows \cc: g_*$ and the fact that $g^*,g_*$ are left exact one easily proves that $g^* \circ -: \Lex(\cg_{\omega}^{\op},\cd) \rightleftarrows \Lex(\cg_{\omega}^{\op},\cc): g_*\circ -$ is also an adjunction. Consequently, \eqref{comm1} indeed describes $(1_\cg \boxtimes g)^*$, as desired. 
	
	Observe that, as $g^*$ is left exact and finite limits are computed pointwise both in $\Lex(\cg_\omega^{\op},\cc)$ and in $\Lex(\cg_\omega^{\op},\cd)$, we have that $g^*\circ - $ is left exact. Consequently, $(1_\cg \boxtimes g)^*$ is left exact. 
\end{proof}
\end{prop}

\begin{prop}\label{leftexactalphaembedding}
	Let $\cc$ be a Grothendieck category and consider a regular cardinal $\alpha$ such that $\cc$ is locally $\alpha$-presentable. Consider the flat embedding $f: \cc \simeq \Ind_\alpha(\cc_\alpha) \ra \Mod(\cc_\alpha)$. Then, for any Grothendieck category $\cg$, we have that $(f \boxtimes 1_\cg)^*$ is left exact.
	
	\begin{proof}
		First, observe that
		\begin{equation} \label{indalpha}
			\begin{aligned}
			\Ind_\alpha(\cc_\alpha) \boxtimes \cg &\simeq \Cont(\Ind_\alpha(\cc_\alpha)^\op,\cg)\\ 
			&\simeq \Cocont(\Ind_\alpha(\cc_\alpha),\cg^\op)^\op\\
			&\simeq \Rex_\alpha(\cc_\alpha,\cg^\op)^\op\\
			&\simeq \Lex_\alpha(\cc_\alpha^\op,\cg)
			\end{aligned}
		\end{equation}	
		and similarly, 
		\begin{equation}\label{mod}
		\begin{aligned}
		\Mod(\cc_\alpha) \boxtimes \cg &\simeq \Cont(\Mod(\cc_\alpha)^\op,\cg)\\ 
		&\simeq \Cocont(\Mod(\cc_\alpha),\cg^\op)^\op\\
		&\simeq [\cc_\alpha,\cg^\op]^\op\\
		&\simeq [\cc_\alpha^\op,\cg].
		\end{aligned}
		\end{equation}	
		Now, consider a regular cardinal $\beta$ such that $\cg$ is $\beta$-presentable and hence $\cg \simeq \Ind_\beta(\cg_\beta)$. We show now that the diagram
		\begin{equation}
			\begin{tikzcd}
			{\Cont(\Ind_\alpha(\cc_\alpha)^\op,\cg)} \arrow[rr, "(f \boxtimes 1_\cg)_*" description] \arrow[d, "\phi"' outer sep=4pt,"\simeq" description] &  & {\Cont(\Mod(\cc_\alpha)^\op,\cg)} \arrow[d, "\psi" outer sep=4pt, "\simeq" description] \\
			{\Lex_\alpha(\cc_\alpha^\op,\cg)} \arrow[rr, "i'" description, hook] \arrow[d, "\simeq" description]                &  & {[\cc_\alpha^\op,\cg]} \arrow[d, "\simeq" description]            \\
			{\Lex_{\alpha,\beta}(\cc_\alpha^\op \otimes \cg_\beta^\op, \Mod(k))} \arrow[rr, "i" description, hook]          &  & {\Lex_{-,\beta}(\cc_\alpha^\op \otimes \cg_\beta^\op, \Mod(k))} 
			\end{tikzcd}
		\end{equation}
		is commutative up to isomorphism, where $\Lex_{-,\beta}(\cc_\alpha^\op \otimes \cg_\beta^\op, \Mod(k))$ is the category of $k$-linear functors $\cc_\alpha^\op \otimes \cg_\beta^\op \ra \Mod(k)$ preserving $\beta$-small limits in the second variable, $i$ and $i'$ are the natural embeddings and $\phi$ and $\psi$ are the isomorphisms provided by \eqref{indalpha} and \eqref{mod} respectively. The commutativity of the lower square is immediate. We prove the commutativity of the upper square. Recall from \Cref{remtheothervariable} that the right adjoint $(f \boxtimes 1_\cg)_*$ is given by $- \circ (f^*)^\op$. Then, given $F \in \Cont(\Ind_\alpha(\cc_\alpha)^\op,\cg)$, we have that 
		\begin{equation}
			\psi \circ (f \boxtimes 1_\cg)_* (F) = F \circ (f^*)^\op \circ Y^\op
		\end{equation}
	 	where $Y:\cc_\alpha \ra \Mod(\cc_\alpha)$ is the Yoneda embedding. 
		On the other hand, we have that
		\begin{equation}
		i'\phi(F) = F \circ Y_\alpha^\op
		\end{equation}
		where $Y_\alpha: \cc_\alpha \ra \Ind_\alpha(\cc_\alpha)$ is the corestriction of the Yoneda embedding. As $f_* \circ Y_\alpha = Y$, we have that $f^* \circ Y = f^* \circ f_* \circ Y_\alpha \cong Y_\alpha$. We can thus conclude that $\psi(-\circ (f^*)^\op)(F) \cong i'\phi(F)$, showing the commutativity (up to ismorphism) of the diagram.
		
		Consequently, to show that $(f \boxtimes 1_\cg)^*$ is left exact, it is enough to show that the left adjoint $a: \Lex_{-,\beta}(\cc_\alpha^\op \otimes \cg_\beta^\op, \Mod(k)) \ra \Lex_{\alpha,\beta}(\cc_\alpha^\op \boxtimes_k \cg_\beta^\op, \Mod(k))$ of the functor $i$ is left exact. We have the following commutative diagram
		\begin{equation*}
			\begin{tikzcd}
			{\Lex_{\alpha,\beta}(\cc_\alpha^\op \otimes \cg_\beta^\op,\Mod(k))} \arrow[rr, "i" description, hook] \arrow[rd, "j" description, hook] &                                                                                                                     & {\Lex_{-,\beta}(\cc_\alpha^\op \otimes_k \cg_\beta^\op,\Mod(k))} \arrow[ld, "j'" description, hook] \arrow[ll, "a" description, bend right=13] \\
			& \Mod(\cc_\alpha \otimes_k \cg_\beta) \arrow[ru, "s'" description, bend right=19] \arrow[lu, "s" description, bend left=19] &                                                                                                                                                         
			\end{tikzcd}
		\end{equation*}
		where $j,j'$ are the natural embeddings. Observe that the adjunction $s\dashv j$ is, up to isomorphism, given by $\Mod(\cc_\alpha) \boxtimes \Mod(\cg_\beta) \rightleftarrows \Ind_\alpha(\cc_\alpha) \boxtimes \Ind_\beta(\cg_\beta)$, which is precisely the tensor product $\boxtimes$ of the left exact reflections $\Ind_\alpha(\cc_\alpha) \subseteq \Mod(\cc_\alpha)$ and $\Ind_\beta(\cg_\beta) \subseteq \Mod(\cg_\beta)$ as in  \cite[\S2.6]{lowenramosgonzalezshoikhet18} and hence again a left exact reflection of $\Mod(\cc_\alpha \otimes_k \cg_\beta)$. Consequently, $s$ is left exact. As $a \cong as'j' \cong sj'$ and $s$ and $j'$ are both left exact, we conclude the argument.
	\end{proof}
\end{prop}

\begin{prop}\label{leftexactgeneralcase}
Let $g:\cc \ra \cd$ be a flat morphism and $\cg$ a Grothendieck category. Let $\alpha$ be a regular cardinal such that $\cg$ is locally $\alpha$-presentable. Denote by $f: \cg \simeq \Ind_\alpha(\cg_\alpha) \ra \Mod(\cg_{\alpha})$ the flat embedding. Then, the left adjoint $(1_\cg \boxtimes g)^*$ is given, up to isomorphism, by the following composition:
\begin{equation}\label{decompositionleftadjoint}
	\begin{tikzcd}
	\Mod(\cg_\alpha) \boxtimes \cc   \arrow[dd, "(f \boxtimes 1_{\cc})^*" description] &  &  & \Mod(\cg_\alpha)  \boxtimes \cd \arrow[lll, "(1_{\Mod(\cg_\alpha)} \boxtimes g)^*" description]                                                             \\
	&&&\\
	  \Ind_\alpha(\cg_\alpha) \boxtimes \cc                                     &   & &  \Ind_\alpha(\cg_\alpha)  \boxtimes \cd \arrow[lll, "(1_\cg \boxtimes g)^*" description, dashed] \arrow[uu, "(f \boxtimes 1_\cd)_*" description]
	\end{tikzcd}
\end{equation}
In addition, $(1_\cg \boxtimes g)^*$ is left exact.
\begin{proof}
	Observe that the following diagram given by the right adjoints
	\begin{equation*}
		\begin{tikzcd}
			\Mod(\cg_\alpha) \boxtimes \cc  \arrow[rrr, "g_*\circ - " description] &  &  & \Mod(\cg_\alpha)  \boxtimes \cd  \\
			&  &  &                                                                            \\
			\Ind_\alpha(\cg_\alpha) \boxtimes \cc \arrow[rrr, "g_*\circ -" description]      \arrow[uu, " -\circ (f^*)^\op" description]                                 &  &  & \Ind_\alpha(\cg_\alpha)  \boxtimes \cd   \arrow[uu, "- \circ (f^*)^\op" description]                                  
		\end{tikzcd}                                                            
	\end{equation*}
	is commutative and hence the diagram
	\begin{equation*}
		\begin{tikzcd}
			\Mod(\cg_\alpha) \boxtimes \cc    \arrow[dd, "(f \boxtimes 1_{\cc})^*" description]                                                       &  &  & \Mod(\cg_\alpha)  \boxtimes \cd  \arrow[lll, "(1_{\Mod(\cg_\alpha)} \boxtimes g)^*" description]    \arrow[dd, "(f \boxtimes 1_\cd)^*" description]                                                    \\
			&  &  &                                                                                                                                          \\
			\Ind_\alpha(\cg_\alpha) \boxtimes \cc   &  &  & \Ind_\alpha(\cg_\alpha)  \boxtimes \cd   \arrow[lll, "(1_\cg \boxtimes g)^*" description] 
		\end{tikzcd}
	\end{equation*}
	given by the left adjoints commutes up to isomorphism. Moreover, as shown in the proof of \Cref{leftexactalphaembedding}, $(f \boxtimes 1_\cd)_* = - \circ (f^*)^\op$ is fully faithful, and hence
	\begin{equation*}
		\begin{aligned}
			(1_\cg \boxtimes g)^* &\simeq (1_\cg \boxtimes g)^* \circ (f \boxtimes 1_\cd)^* \circ (f \boxtimes 1_\cd)_* \\ 
			&\simeq (f \boxtimes 1_\cc)^* \circ (1_{\Mod(\cg_\alpha)} \boxtimes g)^* \circ (f \boxtimes 1_\cd)_*
		\end{aligned}
	\end{equation*} 
	proving the commutativity up to isomorphism of \eqref{decompositionleftadjoint}.
	
	It remains to prove that $(1_\cg \boxtimes g)^*$ is left exact. By \Cref{leftexactfinitelypresentable} we have that $(1_{\Mod(\cg_\alpha)} \boxtimes g)^*$ is left exact, and so is $(f \boxtimes 1_\cd)_*$ because it is a right adjoint. In addition, as a consequence of \Cref{leftexactalphaembedding}, we have that $(f \boxtimes 1_\cd)^*$ is also left exact. Consequently, $(1_\cg \boxtimes g)^*$ is left exact as we wanted to show.
\end{proof}
\end{prop}
This result shows that $\grothflat$ is indeed a monoidal structure with the tensor product of Grothendieck categories. We have the following.
\begin{thm}
	The tensor product $\boxtimes$ of Grothendieck categories induces a monoidal structure 
	$$\boxtimes:\grothflat \times \grothflat \ra \grothflat$$
	such that, given flat morphisms $f:\ca \ra \cb$ and $g:\cc \ra \cd$, the left adjoint $(f \boxtimes g)^*: \ca \boxtimes \cc \ra \cb \boxtimes \cd$ coincides with $f^* \boxtimes g^*$ seen as morphisms in $\groth$.
	\begin{proof}
		The result follows from \S\ref{descrightadjoint} together with \Cref{leftexactgeneralcase}.
	\end{proof}
\end{thm}

The following property of the monoidal structure in $\grothflat$ will be useful further on. 
\begin{prop}\label{preserveembeddings}
	Let $\ca,\cb, \cc$ be Grothendieck categories and $\iota: \cb \ra \cc$ a flat embedding. Then $\iota \boxtimes \ca: \cb \boxtimes \ca \ra \cc \boxtimes \ca$ is a flat embedding.
	\begin{proof}
		We know that $\iota \boxtimes \ca$ is given by $\iota_* \circ -: \Cont(\ca^\op,\cb) \ra \Cont(\ca^\op,\cc)$. One then concludes by the fact that $\iota_*$ is fully faithful.
	\end{proof}
\end{prop}

\subsection{Exponentiable Grothendieck categories}\label{subsecexponentiable}
The tensor product of Gro\-then\-dieck categories $\boxtimes$ endows $\grothflat$ with a monoidal structure. We recall the definition of exponentiable object for this particular monoidal structure.
\begin{defn}
	A Grothendieck category $\ce$ is called \emph{exponentiable} if, for every Grothendieck category $\cb$, there exists a Grothendieck category $\cb^{\ce}$ such that 
	\begin{equation}
	\grothflat(\ca \boxtimes \ce,\cb) \simeq \grothflat(\ca,\cb^{\ce})
	\end{equation}
	for any Grothendieck category $\ca$, and such equivalence is pseudo-natural in $\ca$.
\end{defn}
The rest of the paper will be devoted to providing a characterization of the exponentiable objects in $(\grothflat,\boxtimes)$. But before that, we can already provide an important tool in order to check exponentiability of Grothendieck categories, that exemplifies very well some of the strategies that we will adopt in the rest of the paper and that will be essential in the proof of our main theorem.
\begin{lem}\label{lem:reducingthebase}
	A Grothendieck category $\ce$ is exponentiable if an only if, for every small linear category $\AAA$, the exponential $\Mod(\hat \AAA)^\ce$ exists, where $\hat{\AAA}$ denotes the free completion of $\AAA$ under finite limits. 
	\begin{proof}
		One implication is trivial, thus we only discuss the other. Assume that $\Mod(\hat \AAA)^\ce$ exists for any small $\AAA$ and let us show that $\cd^\ce$ exists for any Grothendieck category $\cd$. We claim that any such $\cd$ can be written as a pullback of a diagram of categories of linear presheaves as follows 
		\begin{center}
			\begin{tikzcd}
				\cd \arrow[dd, "i" description, dashed] \arrow[rr, dashed] &  & \Mod(\hat \CCC) \arrow[dd, "g" description] \\
				&  &                                              \\
				\Mod(\hat{\AAA}) \arrow[rr, "f" description]                        &  & \Mod(\hat \BBB).                            
			\end{tikzcd}
		\end{center}
		Indeed, using the Gabriel-Popescu theorem, given a small subcategory $\AAA$ of $\cd$ generating $\cd$, there exists a flat morphism $\mathsf{Sh}: \Mod(\hat \AAA) \rightleftarrows \cd : i$, which can be specified by inverting a set of arrows $S$ in ${\Mod(\hat \AAA)}$. Define $\{\ra\}$ as the linear category with two objects $a,b$ and morphisms $\{\ra\}(a,a) = \{\ra\}(b,b) = \{\ra\}(a,b) = k$, $\{\ra\}(b,a) = 0$. Furthermore, define $\{\rightleftarrows\}$ as the linear category with two objects $a,b$ such that $\{\rightleftarrows\}(a,a) = \{\rightleftarrows\}(b,b) = \{\rightleftarrows\}(a,b) = \{\rightleftarrows\}(b,a) = k$. Now call $\CCC$ the category $\coprod_S \{\ra\}$ and $\BBB$ the category $\coprod_S \{\rightleftarrows\}$. Then, it is easy to see that one can arrange these data such that the diagram above witnesses $\cd$ as the pullback in $\grothflat$. Now, in order to show that $\cd^\ce$ exists, it is enough to observe that exponentiation commutes with pullbacks and thus we can pointwise compute the exponential over the categories of linear presheaves, and then take the pullback of the resulting diagram, which exists because $\grothflat$ is closed under pullbacks (see \cite[Thm 3.6]{pitts85} for the nonlinear parallel of this result). This pullback will hence satisfy the universal property of $\cd^\ce$ by construction.
	\end{proof}
\end{lem}
In particular, this result allows us to already provide an important class of exponentiable Grothendieck categories. 
\begin{prop}\label{propmodulesareexponentiable}
	Given a small linear category $\DDD$, we have that $\Mod(\DDD)$ is exponentiable.

\end{prop}

\begin{proof}
	By \Cref{lem:reducingthebase}, it is enough to show that for every small linear category $\AAA$, the exponential $\Mod(\hat \AAA)^{\Mod(\DDD)}$ exists. We claim that \[\Mod(\hat \AAA)^{\Mod(\DDD)} \simeq \Mod(\widehat{\AAA \otimes\DDD^\op}).\]
	To show this, we prove that $\Mod(\widehat{\AAA \otimes\DDD^\op})$ satisfies the desired universal property:
	\begin{equation}
		\begin{aligned}
		\grothflat(\cc \boxtimes \Mod(\DDD), \Mod(\hat \AAA)) &\simeq \grothflat([\DDD^\op,\cc],\Mod(\hat \AAA))\\
		&\simeq \Cocont\Lex(\Mod(\hat \AAA),[\DDD^\op,\cc])\\
		&\simeq [\AAA \otimes \DDD^\op, \cc]\\
		&\simeq \Cocont\Lex(\Mod(\widehat{\AAA \otimes\DDD^\op}),\cc)\\
		&\simeq\grothflat(\cc, \Mod(\widehat{\AAA \otimes\DDD^\op})).
		\end{aligned}
	\end{equation}

\end{proof}

\section{$\Mod(k)[\obj]$}

In this brief section we study the properties of the forgetful functor \[\mathbb{U}: \grothflat^{\op} \ra \Cat_k\] and we show that it is representable. Notice that $\mathbb{U}(f) = f^*$. Our treatment is informed of the classical topos theory, where the analogous forgetful functor $\text{Topoi}^\op \to \Cat$ has been studied and shown to be representable \cite[B.3.2.9]{elephant1}. In that case, the functor is represented by the copresheaf topos over finite sets $\Set[\obj]$ and a very similar intuition is fruitful also for our purposes. In the language of topos theory, since a geometric morphism $\ce \to \Set[\obj]$ corresponds precisely to an object of $\ce$, such a topos classifies the geometric \textit{theory of objects}. In analogy with this intuition, the Grothendieck category that we find will be denoted by $\Mod(k)[\obj]$.

\begin{notat}
	Given a module category $\Mod(\AAA)$ with $\AAA$ a small linear category, the subcategory of finitely presentable objects $\Mod(\AAA)_\omega$ of $\Mod(\AAA)$ is usually denoted in the literature by $\fpmod(\AAA)$. We will follow this convention throughout the rest of the paper.
\end{notat}

\begin{prop}\label{forgetfulrep}
	The forgetful functor $\mathbb{U}: \grothflat^{\op} \ra \Cat_k$ is represented by:
	\begin{equation}
		\Mod(k)[\obj] = \Mod(\fpmod(k)^{\op}).
	\end{equation}
	\begin{proof}
		Given any Grothendieck category $\ca$, we have that
		\begin{equation}\label{eqobjectclassifier1}
			\begin{aligned}
			\grothflat(\ca,\Mod(\fpmod(k)^{\op})) &\simeq \Cocont\Lex(\Mod(\fpmod(k)^{\op}),\ca)\\
			&\simeq \Lex [\fpmod(k)^{\op},\ca],
			\end{aligned}
		\end{equation}
		where the second equivalence follows from the fact that $\Mod(\fpmod(k)^{\op})$ is the free cocompletion of $\fpmod(k)^{\op}$ \cite[Thm 4.50]{kelly82basic} and the fact that the left Kan extension of a left exact functor is left exact \cite[Thm 6.12]{kelly82structures}. Now observe that
		\begin{equation}\label{eqobjectclassifier2}
			\Lex [\fpmod(k)^{\op},\ca] \simeq [\{k\},\ca] \simeq \ca.
		\end{equation}
		where $\{k\}$ denotes the full $k$-linear subcategory of $\fpmod(k)^{\op}$ with object the $k$-module $k$. This is as a direct consequence of the fact that $\fpmod(k)$ is the completion under finite colimits of $\{k\}$ in $\Mod(k)$ and hence we can apply \cite[Thm 5.35]{kelly82basic}. Combining \eqref{eqobjectclassifier1} and \eqref{eqobjectclassifier2}, we conclude that the forgetful functor $\grothflat^{\op} \ra \Cat_k$ is represented by $\Mod(\fpmod(k)^{\op})$ as we wanted to show.
	\end{proof}
\end{prop}

\section{Continuous linear categories and injective Grothendieck categories}

 This section studies injective Grothendieck categories, continuous linear categories and then connects the two concepts. These will be relevant technical tools for our main theorem. The relevance of continuous categories in connection to exponentiabilty was pointed out in \cite{johnstonejoyal82}, strongly inspired by the seminal works \cite{hyland81} and \cite{daykelly70}. The approach of Johnstone and Joyal has later inspired \cite{anellejay18} and has influenced our treatment quite deeply.

 The subsection below introduces the notion of (quasi-)injective Grothendieck category. Then we discuss the notion of continuous linear category, finally in the last subsection we see how continuous linear categories and (quasi-)injective Grothendieck categories form bi-equivalent $2$-categories.

\subsection{(Quasi-)Injective Grothendieck categories} \label{qinjgcats}

\begin{defn} An object $\ca$ in $\grothflat$ is \emph{injective} if, given a span like the one in the diagram below

\begin{center}
\begin{tikzcd}
\ca                                                         &  &                                           \\
                                                            &  &                                           \\
\cb \arrow[uu, "f" description] \arrow[rr, "g" description] &  & \cc \arrow[lluu, "h" description, dashed]
\end{tikzcd}
\end{center}
where $g$ is a flat embedding (i.e. its direct image is fully faithful), there exists a dotted flat morphism  $h: \cc \to \ca$ and makes the diagram commute.
\end{defn}

\begin{rem}
Instances of the concept of injectivity are studied everywhere in mathematics, very often with respect to the relevant concept of monomorphism, as in this  case.
\end{rem}

\begin{prop} \label{prop:Modinjective}
 When $\AAA$ is a category with finite limits, $\Mod(\AAA)$ is injective.
\end{prop}
\begin{proof}
According to the definition, and the notations in the diagram below, we need to define a flat morphism $h: \cc \to \Mod(\AAA)$. 
\begin{center}
\begin{tikzcd}
\Mod(\AAA)                                                         &  &                                           \\
                                                            &  &                                           \\
\cb \arrow[uu, "f" description] \arrow[rr, "g" description] &  & \cc \arrow[lluu, "h" description, dashed]
\end{tikzcd}
\end{center}
We will define its left adjoint $h^*: \Mod(\AAA) \to \cc$, as described in the the diagram below, by the formula $h^*:=\Lan_{Y_{\AAA}}(g_*f^*Y_\AAA).$

\begin{center}
\begin{tikzcd}
\AAA \arrow[rd, "Y_\AAA" description] \arrow[rrrddd, "g_*f^*Y_\AAA" description, bend left=49] &                                                                                                       &  &     \\ 
& \Mod(\AAA) \arrow[dd, "f^*" description] \arrow[rrdd, "\Lan_{Y_{\AAA}}(g_*f^*Y_\AAA)" description, dashed] &  &     \\
&   &  &     \\
& \cb \arrow[rr,"g_*"description]                                                                     &  & \cc
\end{tikzcd}
\end{center}
Indeed $h^*$ is cocontinuous (by the universal property of the presheaf construction) and left exact (because $g_*f^*Y_{\AAA}$ is a composition of left exact functors). In order to finish the proof we need to show that (passing to the inverse images), \[f^*\cong g^*h^*.\] Let's check this directly. Recall that, since $g$ is a flat embedding, we have that $g^*g_* \cong 1$.

\[g^*h^*= g^*\Lan_{Y_{\AAA}}(g_*f^*Y_{\AAA}) \cong \Lan_{Y_{\AAA}}(g^*g_*f^*Y_{\AAA}) \cong \Lan_{Y_{\AAA}}(f^*Y_\AAA) \cong f^*. \]

\end{proof}

\begin{prop}
The injective Grothendieck categories are precisely the retracts of the categories of linear presheaves over categories with finite limits.
\end{prop}
\begin{proof}
Assume that $\ca$ is injective. By Gabriel-Popescu, $\ca$ has a flat embedding into a presheaf category $\ca \to \Mod(\AAA)$. We can assume that $\AAA$ has finite limits, indeed Gabriel-Popescu theorem is based on the fact that $\AAA$ is a generating family in $\ca$, and any generating family can be closed under finite limits in $\ca$. 
\begin{center}
\begin{tikzcd}
\ca                                             &  &                                 \\
                                                &  &                                 \\
\ca \arrow[uu, "1"] \arrow[rr, "i" description] &  & \Mod(\AAA) \arrow[lluu, dashed]
\end{tikzcd}
\end{center}
By its injectivity property, this embedding must split. The other implication follows directly from \Cref{prop:Modinjective} and the fact that a rectract of an injective Grothendieck category is injective.
\end{proof}

The notion of quasi-injective Grothendieck category is inspired by the characterization above. 

\begin{defn} 
An object $\ca$ in $\grothflat$ is \emph{quasi-injective} if it is a retract of a category of linear presheaves. The full subcategory spanned by quasi-injective Grothendieck categories will be indicated by \QInj.
\end{defn}

\begin{rem}\label{ccqinj}
Notice that $\QInj$ is biequivalent to the (pseudo)Cauchy-completion\footnote{Also known as Karoubi completion, or completion under pseudo-idempotents.} of $\Presh$, the $2$-category of categories of linear presheaves with flat morphisms. This follows directly from the following proposition.
\end{rem}

\begin{prop}
Let $\ck$ be a Cauchy complete category and $\ca \subset \ck$ be a full subcategory. Then the Cauchy completion $\bar{\ca}$ of $\ca$ coincides with the closure of $\ca$ under idempotents in $\ck$.
\end{prop}
\begin{proof}
In a nutshell, this follows by the fact that idempotents are absolute colimits. Call $\hat \ca$ the closure of $\ca$ in $\ck$ under idempotents. We want to show that $\bar \ca \simeq \hat \ca$. Indeed, by the fact that $\hat \ca$ splits idempotents, we have a functor $\bar \ca \to \hat \ca$ which is given by the universal property of $\bar \ca$. This functor is clearly fully faithful. To show that it is essentially surjective, one uses the explicit presentation of an idempotent. 
\end{proof}

\subsection{Continuous linear categories} \label{contcat}

The theory of continuous categories was initiated in poset theory, where continuous lattices emerged naturally as those complete lattices for which taking suprema of directed subsets commutes with taking infima of arbitrary subsets \cite{zbMATH03379785}. They became relevant in the study of exponentiable objects in the category of topological spaces and in the category of locales after the seminal works of Day and Kelly \cite{daykelly70} and Hyland \cite{hyland81}. A topological space is exponentiable if and only if its frame of opens is a continuous lattice, and more generally a locale is exponentiable if and only if it is continuous. The same intuition has led to the definition of continuous topos \cite{johnstonejoyal82} and continuous $\infty$-topos \cite{anellejay18}.

\begin{defn}\label{def:continuouscat} 
Let $\ca$ be a linear category with filtered colimits, or equivalently, let $\ca$ be a linear category such that the canonical functor $i: \ca \hookrightarrow \Ind(\ca)$ has a left adjoint $\colim: \Ind(\ca) \ra \ca$ (see the first two pages of \cite{zbMATH00825324} for a conceptual discussion about this result). We say that $\ca$ is \emph{continuous} if the left adjoint $\text{colim}: \Ind(\ca) \ra \ca$ has a further left adjoint $l:\ca \ra \Ind(\ca)$, that is, if we have an adjoint triple
\begin{center}
\begin{tikzcd}
\ca \arrow[rr, "i" description, shift right=4] \arrow[rr, "l" description, shift left=5] &  & \text{Ind}(\ca). \arrow[ll, "\text{colim}" description]
\end{tikzcd}
\end{center}
In particular, as $i$ is fully faithful, so is $l$.
\end{defn}

\begin{rem}\label{rem:colimitfunctor}
As hinted by the notation, the functor $\colim$ consists of sending ind-objects, which are formal filtered colimits, to the corresponding actual filtered colimits in $\ca$. The existence of a left adjoint $l$ is an evidence of a good (exceptional) interplay between limits and filtered colimits.
\end{rem}

\begin{exa}\label{exa:finitelypresentable-continuous}
	The most natural example of continuous linear categories are finitely accessible linear categories. To show this, we follow the argument of \cite[Prop 2.4]{johnstonejoyal82}.

\end{exa}
	\begin{proof}
	Let $\ca$ be a finitely accessible linear category. In particular, we have that $\ca$ has filtered colimits and $\ca \simeq \Ind(\ca_\omega)$. Denote by $\iota: \ca_\omega \hookrightarrow \Ind(\ca_\omega)$ the canonical embedding. We claim that $\Ind(\iota): \Ind(\ca_\omega) \ra \Ind(\Ind(\ca_\omega))$ is the left adjoint of the functor $\colim: \Ind(\Ind(\ca_\omega)) \ra \Ind(\ca_\omega)$. First, it is easy to see that $\colim \Ind(\iota) \cong 1_{\Ind(\ca_\omega)}$; this follows from the definition of $\Ind(\iota)$ and the fact that every object in $\Ind(\ca_\omega)$ can be (formally) written as a colimit of objects $A \in \ca_\omega$. This will be the unit of the adjuntion. We now proceed to build the counit $\Ind(\iota) \colim \ra 1_{\Ind(\Ind(\ca_\omega))}$. 
	Consider an ind-object $ X = \mlqq\underset{i \in I}{\colim}\mrqq X_i \in \Ind(\Ind(\ca_\omega))$ where $X_i \in \Ind(\ca_\omega)$. If we apply $\colim$ we obtain the element $\colim_{i \in I} X_i \in \Ind(\ca_\omega)$. In particular, as this is an element of $\Ind(\ca_\omega)$, it can be canonically written as an ind-object $\mlq\underset{j \in J}{\colim}\mrq A_j$ where $A_j \in \ca_\omega$. As the elements of $\ca_\omega$ are finitely presented in $\Ind(\ca_\omega)$, we have that the canonical morphisms
	$$A_j \ra \colim_{i \in I} X_i$$
	factor through some $X_{i_j}$ for some $i_j \in I$. Consequently, we have a morphism $$\mlqq\underset{j \in J}{\colim}\mrqq A_j \ra \mlqq\underset{i \in I}{\colim}\mrqq X_i$$
	in $\Ind(\Ind(\ca_\omega))$ which becomes the identity after applying $\colim$. Observe that 
	\begin{equation*}
		\Ind(\iota)\colim(\mlqq\underset{i \in I}{\colim}\mrqq X_i) = \mlqq\underset{j \in J}{\colim}\mrqq A_j
	\end{equation*}
	and that the morphism provided is the component of a natural transformation $\Ind(\iota) \colim \Rightarrow 1_{\Ind(\Ind(\ca_\omega))}$. The fact that the unit and counit satisfy the triangular equivalences is left to the reader. 
	\end{proof}

\begin{prop} \label{prop:FirstCauchycompPresh}
	A linear category is continuous if and only if it is a retract of a finitely (class\footnote{The notion of class-accessible category will not be introduced in this paper, as it is very intuitive, it is exactly the same of the notion of accessibility, but it is allowed to have a proper class of finitely presentable objects. See \cite{CHORNY20122113} for more details.}-)accessible linear category by functors preserving filtered colimits.
\end{prop}
\begin{proof}
	One implication follows from the definition of continuous category  (\Cref{def:continuouscat}). We proceed to prove the other implication following the arguments of \cite[Prop 2.7]{johnstonejoyal82}.  
	Consider a filtered colimit preserving retract $r: \Ind(\cb) \ra\ca$, with filtered colimit preserving section $s:  \ca \rightarrow \Ind(\cb)$. The category $\Ind(\cb)$ has filtered colimits and hence so does $\ca$, as the retract $r$ preserves filtered colimits. Consequently, we have that the canonical embeddings $\iota_{\ca}: \ca \hookrightarrow \Ind(\ca)$ and $\iota_{\Ind(\cb)}: \Ind(\cb) \ra \Ind(\Ind(\cb))$ have left adjoints $\colim_{\ca}$ and $\colim_{\Ind(\cb)}$ respectively. From the fact that $\colim_{\ca}, \colim_{\Ind(\cb)}$, are colimit preserving and the fact that $r$ and $s$ preserve filtered colimits, together with the definition of $\Ind(r)$ and $\Ind(s)$, one can easily conclude that the following diagram
	\begin{equation*}
		\begin{tikzcd}
			\ca \arrow[rr, "s" description]                                                      &  & \Ind(\cb) \arrow[rr, "r" description]                           &  & \ca                                                \\
			&  &                                                                                                                   &  &                                                    \\
			\Ind(\ca) \arrow[rr, "\Ind(s)" description] \arrow[uu, "\colim_{\ca}" description] &  & \Ind(\Ind(\cb)) \arrow[rr, "\Ind(r)" description] \arrow[uu, "\colim_{\Ind(\cb)}" description] &  & \Ind(\ca) \arrow[uu, "\colim_{\ca}" description]
		\end{tikzcd}
	\end{equation*}
	is commutative. We know by \Cref{exa:finitelypresentable-continuous} that $\Ind(\cb)$ is continuous and hence the functor $\colim_{\Ind(\cb)}$ has a left adjoint $l: \Ind(\cb) \ra \Ind(\Ind(\cb))$. Then, applying (the linear parallel) to \cite[Lem 2.6]{johnstonejoyal82}, one concludes that $\colim_{\ca}$ admits a left adjoint, proving the continuity of $\ca$.  Roughly, the argument goes as follows. A natural candidate for the left adjoint is given by the composition $a \coloneqq \Ind(r) l s$, with the unit $\mu: 1_\ca \Rightarrow \colim_\ca a$ and counit $\epsilon: a \colim_{\ca} \Rightarrow 1_{\Ind(\ca)}$ induced from those of the adjunction $l \dashv \colim_{\Ind(\cb)}$. However, while the composition $(\colim_{\ca}\epsilon) \circ (\mu \colim_{\ca})$ is the identity on $\colim_{\ca}$, we have that the composition $(\epsilon a) \circ (a \mu)$ is not the identity on $a$. Nonetheles, by an argument of Paré (see \cite[\S IV.1, Ex 4]{maclane71} or \cite[Cor 8]{street96}), these data guarantee the existence of a left adjoint of $\colim_{\ca}$. In order to obtain the left adjoint, one can check that though $(\epsilon a) \circ (a \mu)$ is not the identity, it is idempotent and, as idempotents are split in the category of linear functors $[\ca, \Ind(\ca)]$ (because they are in $\Ind(\ca)$), we obtain an splitting $a \overset{\sigma}{\Rightarrow} l' \overset{\tau}{\Rightarrow} a$. One can then check that $l': \ca \ra \Ind(\ca)$ is the left adjoint of $\colim_{\ca}$; more concretely, the unit and counit of the adjunction are respectively given by $(\colim_{\ca}\sigma) \circ \mu$ and $\epsilon \circ (\tau \colim_{\ca})$.
\end{proof}

\begin{rem}
	Observe that given a continuous linear category $\ca$ which is locally small, but not small, we can only a priori write it as a retract of a category $\Ind(\cb)$ where $\cb$ is also locally small, but not small in general. We will be interested in the continuous categories $\ca$ that are determined by small data, or in other words, that can be written as a retract of a finitely accessible linear category $\Ind(\AAA)$ where $\AAA$ is small. For this purpose we introduce in what follows the notion of \emph{standard presentation} of a continuous category. 
\end{rem}

\begin{notat}
Let $\ca$ be a continuous linear category and $\AAA \subseteq \ca$ a full subcategory. In particular, we call $\epsilon: \Ind(\AAA) \ra \ca$ the functor given by the composite
\begin{equation*}
	\begin{tikzcd}
		\Ind(\AAA) \arrow[r, hook]  &\Ind(\ca) \arrow[r, "\colim"] & \ca.
	\end{tikzcd}
\end{equation*}
Notice that this functor preserves filtered colimits because both components do. The following result provides a sufficient condition for this functor to have a left adjoint $\beta: \ca \ra  \Ind(\AAA)$.
\end{notat} 

\begin{prop}\label{prop:sufficientconditionstandard}
	Let $\ca$ be a continuous linear category and $\AAA \subseteq \ca$ a full subcategory. Assume that there exists a functor $h: \ca \ra \Ind(\AAA)$ such that the composition
	\begin{equation}
			\begin{tikzcd}
				\ca \arrow[r, "h"] & \Ind(\AAA) \arrow[r, hook] & \Ind(\ca) \arrow[r, "\colim"] & \ca
			\end{tikzcd}
	\end{equation}
	is isomorphic to the identity, or in other words, that every object in $\ca$ can be expressed as a filtered colimits of elements in $\AAA$. Then, the filtered colimit preserving functor $\epsilon: \Ind(\AAA) \ra \ca$ as defined above has a fully faithful left adjoint.
	\begin{proof}
	The argument we will follow can be extracted from the proof of (i) $\Rightarrow$ (iii) in \cite[Prop C.4.2.18]{johnstonejoyal82}.
	
	Consider the composition 
	\begin{equation*}
		\begin{tikzcd}
			\ca \arrow[r, "l"] & \Ind(\ca) \arrow[r, "\Ind(h)"] & \Ind(\Ind(\AAA)) \arrow[r, "\colim"] & \Ind(\AAA) \arrow[r,hook] &\Ind(\ca)
		\end{tikzcd}
	\end{equation*}	
	and denote it by $k: \ca \ra \Ind(\ca)$. Given $A \in \ca$, we are going to show that $k(A) \cong l(A)$. If $l(A)$ is the ind-object $\mlq\underset{i \in I}{\colim}\mrq A_i$ in $\Ind(\ca)$, and $h(A_i)$ is the ind-object $\mlq\underset{j \in J_i}{\colim}\mrq B_{ij}$ in $\Ind(\AAA)$, we have that $k(A)$ is given by the ind-object $\underset{(i,j) \in I\times J_i}{\mlq\colim\mrq} B_{ij}$ obtained when computing $\colim_{i \in I}(\mlq\underset{j \in J_i}{\colim}\mrq B_{ij})$ in $\Ind(\ca)$. Observe that after applying the functor $\colim: \Ind(\ca) \ra \ca$ to the natural morphisms $B_{ij} \ra h(A_i)$ in $\Ind(\ca)$ we obtain morphisms 
	\begin{equation*}
		 B_{ij} \ra A_i
	\end{equation*}
	in $\ca$ which induce a morphism between ind-objects
	\begin{equation*}
		k(A) = \underset{(i,j) \in I \times J_i}{\mlq\colim\mrq} B_{ij} \xrightarrow{f} \mlq\underset{i \in I}{\colim}\mrq A_i = l(A)
	\end{equation*} 
	in $\Ind(\ca)$. On the other hand, an easy computation shows that $A \cong \colim(k(A)) \in \ca$, and thus, applying the adjunction $l \dashv \colim$ (which we have because $\ca$ is continuous by hypothesis), we obtain a morphism 
	\begin{equation*}
		l(A) \xrightarrow{g} k(A).
	\end{equation*}
	An easy computation shows that $f$ and $g$ are inverse to each other and hence $k(A) \cong l(A)$. This in particular shows that $l: \ca \ra \Ind(\ca)$ factors, up to isomorphism, through $\Ind(\AAA) \hookrightarrow \Ind(\ca)$. Let $\beta: \ca \ra \Ind(\AAA)$ the functor through which $l$ factors. We therefore have the commutative diagrams	
	\begin{equation*}
		\begin{tikzcd}
			\ca \arrow[d, Rightarrow, no head] &  & \Ind(\ca) \arrow[ll, "\colim" description]                    &  & \ca \arrow[d, Rightarrow, no head] \arrow[rr, "l" description, hook] &  & \Ind(\ca)                  \\
			\ca                       &  & \Ind(\AAA) \arrow[u, hook] \arrow[ll, "\epsilon" description] &  & \ca \arrow[rr, "\beta" description]                         &  & \Ind(\AAA) \arrow[u, hook]
		\end{tikzcd}
	\end{equation*}
	From these diagrams and the adjuntion $l \dashv \colim$ one can readily conclude that $\beta \dashv \epsilon$ is an adjunction and $\beta$ is fully faithful. 
	\end{proof}
\end{prop}

\begin{defn} \label{def:standardpresentation}
Let $\ca$ be a continuous linear category. If there exists a small full linear subcategory $\AAA \subset \ca$ in the hypothesis of \Cref{prop:sufficientconditionstandard}, we call the induced adjunction
\begin{center}
\begin{tikzcd}
\ca \arrow[rr, "\beta" description, hook, shift left=2] &  & \text{Ind}(\AAA) \arrow[ll, "\epsilon" description, shift left=2]
\end{tikzcd}
\end{center}
a \emph{standard presentation} of $\ca$. In particular, $\epsilon$ preserves filtered colimits and $\beta$ is fully faithful.
\end{defn}

A similar notion appeared in \cite{anellejay18}, with the same purpose. Our notion is more flexible and indeed works for non-cocomplete continuous categories.

\begin{rem} \label{rem:standardpresentationfinitelyaccessible}
	Let $\ca$ be a finitely accessible linear category. We know by \Cref{exa:finitelypresentable-continuous} that $\ca$ is continuous. Notice that the small full subcategory $\ca_\omega \subseteq \ca \simeq \Ind(\ca_\omega)$ gives rise to a standard presentation of $\Ind(\ca_\omega)$. Indeed, the functor $h \coloneqq 1_{\Ind(\ca_\omega)}$ trivially satisfies the assumption from \Cref{prop:sufficientconditionstandard}.
\end{rem}

\begin{rem}\label{rem:standardpresentationpresentable}
	Observe that any locally presentable continuous linear category $\ca$ admits a presentation. Indeed, if $\ca$ is locally $\lambda$-presentable, the full dense subcategory $\ca_\lambda \subseteq \ca$ together with the natural morphism
	\begin{equation*}
		h: \ca \ra \Ind(\ca_\lambda): A \mapsto  \underset{(X_i \ra A) \in \ca_\lambda \downarrow A}{\mlq \colim \mrq} X_i  
	\end{equation*}
	satisfy the assumption of \Cref{prop:sufficientconditionstandard} as a direct consequence of \cite[Prop 1.22]{adamek_rosicky_1994}. Notice that this morphism is well-defined because the indexing category $\ca_\lambda \downarrow A$ is $\lambda$-filtered, and thus in particular ($\omega$-)filtered.
\end{rem}
	
\begin{prop} \label{prop:CauchycompPresh}
	A linear category is continuous with a standard presentation if and only if it is a retract by functors preserving filtered colimits of a finitely accessible linear category $\Ind(\BBB)$ with $\BBB$ small.
\end{prop}
\begin{proof}
One implication follows from the definition of standard presentation  (\Cref{def:standardpresentation}). We proceed to prove the other implication following the arguments of \cite[Prop C.4.2.18]{elephant2}.  

Assume now that $\ca$ is a retract of a finitely accessible linear category $\Ind(\BBB)$ with $\BBB$ a small linear category. Then, by \Cref{prop:FirstCauchycompPresh}, we can conclude that $\ca$ is continuous. It remains to show that $\ca$ has a standard presentation, i.e. there exists a small subcategory $\AAA$ of $\ca$ and a functor $h: \ca \ra \Ind(\AAA)$ such that the composition
\begin{equation*}
	\begin{tikzcd}
		\ca \arrow[r, "h"] & \Ind(\AAA) \arrow[r, hook] & \Ind(\ca) \arrow[r, "\colim"] & \ca
	\end{tikzcd}
\end{equation*}
is isomorphic to the identity of $\ca$. As pointed out in \Cref{rem:standardpresentationfinitelyaccessible}, we have that $\BBB \subseteq \Ind(\BBB)$ with the identity functor $\Ind(\BBB) \ra \Ind(\BBB)$ provide a standard presentation of $\Ind(\BBB)$. We claim that if $r: \Ind(\BBB) \ra \ca$ is the given filtered colimit preserving retract, the full subcategory of $r(\BBB) \subseteq \ca$ with objects those in the image of $\BBB$ under $r$ together with the functor
\begin{equation*}
	\begin{tikzcd}
		\ca \arrow[r, "s"] & \Ind(\BBB) \arrow[r, "\Ind(r_{|\BBB})"] & \Ind(r(\BBB))
	\end{tikzcd}
\end{equation*}
provide a standard presentation of $\ca$, where $s: \ca \ra \Ind(\BBB)$ is the corresponding filtered colimit preserving section. To show this we have to check that the composition
\begin{equation*}
	\begin{tikzcd}
		\ca \arrow[r, "s"] & \Ind(\BBB) \arrow[r, "\Ind(r_{|\BBB})"] & \Ind(r(\BBB)) \arrow[r,hook] &\Ind(\ca) \arrow[r, "\colim"] &\ca
	\end{tikzcd}
\end{equation*}
is isomorphic to the identity on $\ca$. Given $A \in \ca$, if $s(A)$ is the ind-object $\mlq \underset{i \in I}{\colim}\mrq X_i$ in $\Ind(\BBB)$, we then have that this composition sends $A$ to $\colim_{i \in I} r(X_i)$. On the other hand, using the retraction and the fact that $r$ preserves filtered colimits, we have that
\begin{equation*}
	A = rs(A) = r(\mlq \underset{i \in I}{\colim}\mrq X_i)= \colim_{i \in I}r(X_i),
\end{equation*}
which concludes the argument. 
\end{proof}

\subsection{Continuous categories and quasi-injective Grothendieck categories} \label{contvsqinj}

In this section we draw a connection between continuous linear categories and quasi-injective Grothendieck categories.

\subsubsection{Presheaves and locally finitely accessible categories}

\begin{thm}\label{thm:equivPreshAcc}
There is an equivalence of $2$-categories between the $2$-category $\Presh$ of categories of linear presheaves (with morphisms the flat morphisms) and the $2$-category $\omega \mhyphen\Acc$ of finitely accessible linear categories (with morphisms the filtered colimit preserving linear functors). The equivalence is induced by taking points:
\[\mathsf{S}: \omega\mhyphen\Acc \rightleftarrows \Presh : \mathsf{pt}.\]
\end{thm}
\begin{proof}
Let us describe the functors in both direction. The construction is identical to the one in the PhD thesis of the first author \cite{diliberti2020scott,thcat}.
\begin{itemize}
	\item[$\mathsf{S}$] this functor is defined by $\ca \mapsto \omega\mhyphen \Acc(\ca, \Mod(k)).$ This definition explains also its action on morphisms, that is for a functor preserving filtered colimits $f: \ca \to \cb$, a canonical precomposition $\mathsf{S}(f)^*(g) = g \circ f$. $\mathsf{S}(f)$ is cocontinuous and lex, and thus is the inverse image of a flat morphism \[\mathsf{S}f^*: \mathsf{S}\cb \leftrightarrows \mathsf{S} \ca : \mathsf{S} f_*.\] It is easy to see that $\omega\mhyphen \Acc(\ca, \Mod(k))$ coincides with $\Mod(\ca_\omega^\op)$. Indeed, calling $i: \ca_\omega \to \ca$ the inclusion of finitely presentable objects, the Kan-restriction paradigm yields an equivalence of categories 
	\[-\circ i: \omega\mhyphen \Acc(\ca, \Mod(k)) \rightleftarrows [\ca_\omega, \Mod(k)] : \Lan_i.\] 
	Thus $\mathsf{S}\ca$ is always a presheaf category, as desired.
	\item[$\mathsf{pt}$] acts as expected, mapping $\Mod(\AAA) \mapsto \grothflat(\Mod(k),\Mod(\AAA))$, which is easily seen to be a finitely accessible linear category. 
	\end{itemize}
The fact that this constructions are one the inverse of the other is an application of an enriched version Diaconescu's theorem, we provide a short proof below.

	\begin{equation}
		\begin{aligned}
		\mathsf{pt}\mathsf{S}(\ca) & \simeq \mathsf{Cocontlex}(\Mod(\ca_\omega^\op), \Mod(k)) \\
		& \simeq  \mathsf{Lex}(\ca_\omega^\circ, \Mod(k)) \\
		& \simeq \mathsf{Ind}(\ca_\omega) \\
		& \simeq \ca \\
		\end{aligned}
	\end{equation} 

	\begin{equation}
		\begin{aligned}
		\mathsf{S}\mathsf{pt}(\Mod(\AAA)) & \simeq \mathsf{S}(\mathsf{Cocontlex}(\Mod(\AAA), \Mod(k))) \\
		& \simeq  \mathsf{S}(\mathsf{Lex}(\AAA, \Mod(k))) \\
		& \simeq \mathsf{S}(\mathsf{Ind}(\AAA^\op)) \\
		& \simeq \Mod(\AAA). \\
		\end{aligned}
	\end{equation} 
\end{proof}

\subsubsection{Continuous linear categories and quasi-injective Grothendieck categories} \label{equi}
\begin{thm}\label{equivalencecontinuousinjective}
There is an equivalence of $2$-categories between the $2$-category $\QInj$ of quasi-injective Grothendieck categories with flat morphisms and the $2$-category $\Cnt$ of continuous linear categories with a standard presentation and filtered colimit preserving linear functors. The equivalence is induced by taking points:
\[\mathsf{S}: \Cnt \rightleftarrows \text{\QInj} : \mathsf{pt},\]
 where $\QInj$  is the subcategory of quasi-injective topoi, those that are retracts of presheaf topoi.
\end{thm}
\begin{proof}
Since $\omega\mhyphen \Acc$ is biequivalent to $\Presh$ by \Cref{thm:equivPreshAcc}, their Cauchy com\-pletions are biequivalent too. We conclude by observing that $\Cnt$ is the Cauchy completion of $\omega\mhyphen \Acc$ by \Cref{prop:CauchycompPresh} and $\QInj$ is the Cauchy completion of $\Presh$ by \Cref{ccqinj}.
\end{proof}

\section{Main theorem}

The main result of our paper is that a Grothendieck category is exponentiable if and only if it is a continuous linear category.

\begin{thm}\label{maintheorem}
A Grothendieck category is exponentiable in $\grothflat$ if an only if it is a continuous category. In particular every finitely presentable Grothendieck category is exponentiable.
\end{thm}

We split the discussion in two subsections, and we start from proving the necessity of the condition.

\subsection{Exponentiable Grothendieck categories are continuous}
\begin{thm} \label{exponentiable are continuous}
Exponentiable Grothendieck categories are continuous.
\end{thm}
\begin{proof}
The proof is locally easy, but globally tricky, thus we split it in little steps.	
\begin{itemize}[leftmargin=4em]
	\item[Step 1] We show that if $\ce$ is exponentiable then $\Mod(k)[\mathbb{O}]^\ce$ is injective. Assume that $\ce$ is exponentiable. Consider the diagram below.
	\begin{center}
	\begin{tikzcd}
	{\Mod(k)[\mathbb{O}]^\ce} &  &                          &  & {\Mod(k)[\mathbb{O}]}               &  &                                        \\
                              &  &                          &  &                                         &  &                                        \\
	\ca \arrow[uu] \arrow[rr]     &  & \cb \arrow[lluu, dashed] &  & \ca \boxtimes \ce \arrow[rr] \arrow[uu] &  & \cb \boxtimes \ce \arrow[lluu, dotted]
	\end{tikzcd}
	\end{center}
	The right hand side of the diagram corresponds to the left one via the adjunction between exponentiation and tensoring. In order to show that the dashed arrow exists, it's enough to show that the dotted one exists. The latter exists because the tensor product preserves flat embeddings as shown in \Cref{preserveembeddings} and $\Mod(k)[\mathbb{O}]$ is a presheaf category over a finitely complete small linear category $\fpmod(k)^\op$ and thus injective by \Cref{prop:Modinjective}.
	\item[Step 2] We have that 
	\begin{equation*}
		\begin{aligned}
			\ce &\simeq \grothflat(\ce,\Mod(k)[\obj])\\ 
			&\simeq \grothflat(\Mod(k),\Mod(k)[\obj]^\ce)\\
			&\simeq \mathsf{pt}(\Mod(k)[\mathbb{O}]^\ce),
		\end{aligned}
	\end{equation*}
	where the first equivalence follows by \Cref{forgetfulrep}.
	\item[Step 3] Since $\Mod(k)[\mathbb{O}]^\ce$ is injective, by \Cref{equivalencecontinuousinjective} and Step 2, $\ce$ is continuous.
\end{itemize}
\end{proof}

\subsection{Continuous Grothendieck topoi are exponentiable}

The proof that continuous Grothendieck topoi are exponentiable is involved, and we will try to keep it as clear as possible. For this reason, let us introduce the general strategy that we will follow:

\begin{enumerate}
	\item We show that $\ce$ is exponentiable if and only if $\Mod(k)[\mathbb{O}]^\ce$ exists.
	\item In order to investigate the proper candidate of $\Mod(k)[\mathbb{O}]^\ce$, we show that if it exists, it must coincide with $\mathsf{S}(\ce)$.
	\item Since the case of locally finitely presentable Grothendieck categories is much easier and of big interest, we show that if $\ce$ is locally finitely presentable, then $\mathsf{S}(\ce)$ has the universal property of $\Mod(k)[\mathbb{O}]^\ce$.
	\item We give the proof in the general case.
\end{enumerate}

\begin{lem} \label{limitexp}
$\ce$ is exponentiable if and only if $\Mod(k)[\mathbb{O}]^\ce$ exists.
\end{lem}
\begin{proof}
One implication is trivial, we concentrate on the other. By \Cref{lem:reducingthebase}, it is enough to show that if $\Mod(k)[\mathbb{O}]^\ce$ exists, then $\Mod(\hat \AAA)^\ce$ exists for every small linear category $\AAA$. 
We claim that 
\[\Mod(\hat \AAA)^\ce \simeq (\Mod(k)[\mathbb{O}]^\ce)^{\Mod(\AAA^\op)}.\] 
Observe this formula is meaningful because $\Mod(\AAA)$ is indeed exponentiable, as proven in \Cref{propmodulesareexponentiable}.
We show that $(\Mod(k)[\mathbb{O}]^\ce)^{\Mod(\AAA^\op)}$ has the correct universal property. Indeed, we have that
	\begin{equation}
		\begin{aligned}
			\grothflat(\ca, (\Mod(k)[\mathbb{O}]^\ce)^{\Mod(\AAA^\op)}) & \simeq \grothflat({(\Mod(\AAA^\op) \boxtimes \ca) \boxtimes \ce } , \Mod(k)[\mathbb{O}]) \\
			& \simeq  \Mod(\AAA^\op) \boxtimes (\ca \boxtimes \ce) \\
			& \simeq [\AAA, \ca \boxtimes \ce] \\
			& \simeq \mathsf{Cocontlex}(\Mod(\hat \AAA), \ca \boxtimes \ce) \\
			& \simeq \grothflat(\ca \boxtimes \ce, \Mod(\hat \AAA)),\\
		\end{aligned}
	\end{equation}
where the second equivalence follows from \Cref{forgetfulrep} and the associativity of $\boxtimes$.

\end{proof}

\begin{thm} \label{computationexp}
If $\Mod(k)[\obj]^\ce$ exists, then it must be $\mathsf{S}(\ce)$.
\end{thm}
\begin{proof}
This proof is based on a slick double counting over $\mathsf{S}\mathsf{pt}\Mod(k)[\obj]^\ce$. We have that
\[
\Mod(k)[\obj]^\ce \stackrel{(a)}{\simeq}  \mathsf{S}\mathsf{pt}\Mod(k)[\obj]^\ce \stackrel{(b)}{\simeq} \mathsf{S}(\ce),
 \]
where the equivalence (a) is justified by the fact that $\Mod(k)[\obj]^\ce$ is injective, as shown in Step 1 of \Cref{exponentiable are continuous}, and (b) follows from the fact that $\mathsf{pt}(\Mod(k)[\obj]^\ce) \simeq \ce$.
\end{proof}

\begin{thm}\label{locfinpresexponentiable}
If $\ce$ locally finitely presentable then  $\mathsf{S}(\ce)$ satisfies the universal property of the exponential $\Mod(k)[\obj]^\ce$, thus $\ce$ is exponentiable.
\end{thm}
\begin{proof}
	\begin{equation}
		\begin{aligned}
		\grothflat(\ca, \mathsf{S}(\ce) ) & \simeq \mathsf{Cocontlex}(\Mod(\ce_\omega^\op), \ca) \\
		& \simeq  \mathsf{Lex}(\ce_\omega^\circ, \ca) \\
		& \simeq \mathsf{Cocont}(\mathsf{Ind}(\ce_\omega), \ca^\circ)^\circ \\
		& \simeq \ca \boxtimes \mathsf{Ind}(\ce_\omega)\\
		& \simeq \grothflat(\ca \boxtimes \ce, \Mod(k)[\obj] ). \\
		\end{aligned}
	\end{equation}
\end{proof}

We are now ready to provide the proof of the general version of our main theorem. Our proof strategy follows closely the one in \cite{anellejay18}, which is just a formal presentation of \cite{johnstonejoyal82}. This last proof is deeply inspired by Hyland's \cite{hyland81}.

\begin{thm} \label{thm:contimpliesexp}
If $\ce$ is a continuous Grothendieck category then $\Mod(k)[\obj]^\ce$ exists and is $\mathsf{S}(\ce)$.
\end{thm}

\begin{proof}
	We will show that $\mathsf{S}(\ce)$ has the correct universal property:
\begin{itemize}[leftmargin=4em]
	\item[Step 1] Because $\ce$ is locally presentable and continuous, we know by \Cref{rem:standardpresentationpresentable} that we have a standard presentation 
\begin{center}
\begin{tikzcd}
\ce \arrow[rr, "\beta" description, hook, shift left=2] &  & \text{Ind}(\DDD) \arrow[ll, "\epsilon" description, shift left=2]
\end{tikzcd}
\end{center}
for some small linear subcategory $\DDD$. Recall that $\beta$ is fully faithful and $\epsilon$ is filtered colimit preserving. As $\beta$ is a left adjoint, it also preserves filtered colimits. Therefore, we can apply $\mathsf{S}$ to both $\beta$ and $\epsilon$, and using \Cref{equivalencecontinuousinjective}, we obtain 
\[\mathsf{S}(\beta) :\mathsf{S} \ce \leftrightarrows \Mod (\DDD^\op) : \mathsf{S}(\epsilon). \]  
Moreover, it follows by functoriality (or again by \Cref{equivalencecontinuousinjective}) that $\mathsf{S}(\epsilon) \circ \mathsf{S}(\beta)$ is isomorphic to the identity. Finally, applying \Cref{equivalencecontinuousinjective} one last time, we know that passing to the points we get back where we started from, that is $\mathsf{ptS}(\beta) \cong \beta$ and similarly for $\epsilon$. Going back to \Cref{thm:equivPreshAcc}, we even have a precise description of the action of $\mathsf{S}(\beta)^*$, this is \[ \mathsf{S}(\beta)^*(-) = - \circ \beta,\] and analogously for $\mathsf{S}(\epsilon)$. 
\item[Step 2] For all Grothendieck categories $\ca$, we have that $\ca \boxtimes \ce$ is reflective in $\ca \boxtimes \text{Ind}(\DDD)$ via precomposition, as shown below \[(-) \circ \epsilon^\circ : \mathsf{Cont}(\ce^\circ, \ca) \leftrightarrows \mathsf{Cont}(\text{Ind}(\DDD)^\circ, \ca) : (-) \circ \beta^\circ.\] Notice that when $\ca$ is $\Mod(k)$, this returns the original adjunction between $\beta$ and $\epsilon$, as $\Mod(k)$ is the unit of the monoidal structure.
	\item[Step 3] We need to show that for every Grothendieck category $\ca$ one has that 
\[\grothflat(\ca, \mathsf{S} \ce ) \simeq \grothflat(\ca \boxtimes \ce, \Mod(k)[\obj] ) \simeq \ca \boxtimes \ce.\]	Indeed, using Step 1, we have that $\grothflat(\ca, \mathsf{S} \ce )$ is reflective in $\grothflat(\ca, \Mod(\DDD^\op) )$ via precompotion with $\mathsf{S}(\beta)$ and $\mathsf{S}(\epsilon)$, \[\mathsf{S}(\beta) \circ (-): \grothflat(\ca, \mathsf{S} \ce ) \leftrightarrows \grothflat(\ca, \Mod(\DDD^\op) ): \mathsf{S}(\epsilon) \circ (-).\]  Now notice that $\grothflat(\ca, \Mod(\DDD^\op) )$ is precisely $\ca \boxtimes \text{Ind} (\DDD)$. 
\item[Step 4] Putting together the previous considerations, we have the following diagram.
\[\begin{tikzcd}
	{\grothflat(\ca, \mathsf{S}\ce)} && {\grothflat(\ca, \Mod(\DDD^\op))} \\
	\\
	{\ca \boxtimes \ce} && {\ca \boxtimes \text{Ind} (\DDD).}
	\arrow[shift left=2, hook, from=3-1, to=3-3,"(-) \circ \epsilon^\op" description]
	\arrow[shift left=2, hook, from=1-1, to=1-3, "\mathsf{S}(\beta) \circ (-)" description]
	\arrow["\simeq"{description}, from=1-3, to=3-3]
	\arrow[shift left=2, two heads, from=1-3, to=1-1, "\mathsf{S}(\epsilon) \circ (-)" description]
	\arrow[shift left=2, two heads, from=3-3, to=3-1, "(-)\circ \beta^\op" description]
\end{tikzcd}\]
If we show  that $\mathsf{S} (\beta) \circ \mathsf{S}(\epsilon) \circ (-)$ has the same action of $(-) \circ \beta^\circ \circ \epsilon^\circ  $ , then they have the same fixed points, and thus we obtain that $\grothflat(\ca, \mathsf{S} \ce )  \simeq \ca \boxtimes \ce$, as desired.
Now, to finish the proof we pass to the inverse images, which provide us a concrete computation of the action,  \[(\mathsf{S} (\beta) \circ \mathsf{S}(\epsilon) \circ (-))^* =  (-)^* \circ \mathsf{S} (\epsilon)^* \circ \mathsf{S}(\beta)^* = (-)^* \circ \epsilon\circ \beta  ,\] which gets identified to $(-) \circ \beta^\circ \circ \epsilon^\circ $ in the equivalence between $\grothflat(\ca, \Mod(\DDD) )$ and $\ca \boxtimes \text{Ind} (\DDD)$. 
\end{itemize}
\end{proof}

\section{Examples from algebraic geometry}
In this section, we show that the category of quasi-coherent sheaves $\Qcoh(X)$ of a quasi-compact quasi-separated $k$-scheme $X$ is exponentiable (see \Cref{qcohexponentiable}). Fixed such a scheme $X$, the rest of the section is then devoted to obtain a description of the exponential $\Qcoh(Y)^{\Qcoh(X)}$ for a scheme $Y$ and provide concrete examples for particular schemes $X$ and $Y$. We proceed by following the next steps:
\begin{enumerate}
	\item We provide a description of $\Mod(k)[\obj]^{\Qcoh(X)}$ for any quasi-compact quasi-separated $k$-scheme $X$ (\Cref{qcohexponentiable}). We compute it in the particular case in which $X$ is a noetherian $k$-scheme and more concretely when $X= \Spec(A)$ with $A$ a noetherian $k$-algebra (\Cref{exanoeth});
	\item Given any quasi-compact quasi-separated scheme $Y$, we provide a description of $\Qcoh(Y)^{\Qcoh(X)}$ as a pullback of categories of linear presheaves of the form $\Mod\left((\DDD \otimes_\omega \coh(X))^\op\right)$ where $\otimes_\omega$ denotes the tensor product of finitely cocomplete linear categories from \cite{lopezfranco13} (\Cref{expgeneraldescription});
	\item By means of an intermediate result (\Cref{modulecategorybase}), we show that the computation of one of the objects of the pullback diagram gets simplified in the case in which $Y$ is a noetherian scheme (\Cref{expmodulecats}); 
	\item To conclude, we show that \Cref{modulecategorybase} can also be applied to compute $\Mod(A)^{\Qcoh(X)}$ when $A$ is a von Neumann regular $k$-algebra (\Cref{fpabsolutelyflat}) and we provide the computation (\Cref{exabsolutelyflat}).
\end{enumerate}
\begin{notat}
	In order to ligthen the notations, from this point on we will not mention the base commutative ring $k$ unless necessary. Therefore, all schemes and algebras are $k$-schemes and $k$-algebras, and all categories are $k$-linear categories.
\end{notat}
\begin{prop}\label{qcohexponentiable}
	Let $X$ be a quasi-compact quasi-separated scheme, then $\Qcoh(X)$ is exponentiable. In particular, we have that 
	\begin{equation}
		\Mod(k)[\obj]^{\Qcoh(X)} \simeq \Mod(\lfp(\Qcoh(X))^\op),
	\end{equation}
	where $\lfp(\Qcoh(X))$ denotes the full subcategory of $\Qcoh(X)$ given by the locally finitely presentable quasi-coherent sheaves, that is, the quasi-coherent sheaves that locally are a cokernel of a morphism between two finite direct sums of the structure sheaf.
	\begin{proof}
		By \cite[Cor 6.9.12]{grothendieckdieudonne71} we know that $\Qcoh(X)$ is a locally finitely presentable linear category, with the finitely presentable objects given by the locally finitely presentable quasi-coherent sheaves (see \cite[Prop 75]{murfet06}). We conclude by applying \Cref{locfinpresexponentiable}.
	\end{proof}
\end{prop}
\begin{exa}\label{exanoeth}
	Let $X$ be a noetherian scheme. Then we have that 
	$$\Mod(k)[\obj]^{\Qcoh(X)} \simeq \Mod(\coh(X)^\op).$$
	Indeed, as $X$ is noetherian, the finitely presented objects of $\Qcoh(X)$ are precisely the coherent sheaves, that is $\Qcoh(X)_\omega \simeq \lfp(\Qcoh(X)) \simeq \coh(X)$ (see, for example, \cite[Tag 01XZ]{stacksproject}). In particular, given a noetherian commutative algebra $A$, we have that $$\Mod(k)[\obj]^{\Mod(A)} \simeq \Mod(\fpmod(A)^\op),$$
	where $\fpmod(A)$ denotes the category of finitely presentable $A$-modules.
\end{exa}

\begin{rem}[A topological detour]
	\Cref{qcohexponentiable} has a very deep geometric meaning, especially from the topos theoretic point of view, which deserves a bit of space to be expanded. It is known to a general topologist that expontentiability in the category $\mathsf{Top}$ of topological spaces is a quite non-trivial problem. When $X$ is locally compact and Hausdorff, the compact open topology on $\mathsf{Top}(X, -)$ shows that $X$ is exponentiable in $\mathsf{Top}$. In addition, it is known that the frame of opens of a locally compact and Hausdorff space is continuous. This observation is a strong motivation to study continuous frames, which somehow correspond to this very well-behaved notion of topological spaces. It is therefore not a surprise that continuous frames are exponentiable among locales, and that continuous topoi are exponentiable among topoi! It is well known that quasi-compact separated schemes are a linear analog of compact Hausdorff spaces, and thus fall under the big scheme of locally compact Hausdorff. Along these lines, the category of sheaves over any topological manifold is continuous and thus exponentiable. We cannot consider this as a strong evidence that every scheme is exponentiable, but to the authors it seems a good indication that any (category of quasi-coherent sheaves over a) scheme might be exponentiable with respect to this monoidal structure. We therefore pose the following open question.
\end{rem}

\begin{quest}
	Is it true that every scheme is exponentiable?
\end{quest}

Let $X$ be a quasi-compact quasi-separated scheme. From \Cref{qcohexponentiable} we know that $\Qcoh(X)$ is exponentiable, and we have an explicit description of the exponential $\Mod(k)[\obj]^{\Qcoh(X)}$. The rest of the section will be devoted to providing a description of the exponential $\Qcoh(Y)^{\Qcoh(X)}$ for a quasi-compact quasi-separated scheme $Y$ and some examples.
\begin{prop}\label{expgeneraldescription}
	Let $X,Y$ be a quasi-compact quasi-separated schemes. Then $\Qcoh(Y)^{\Qcoh(X)}$ can be expressed as a pullback in $\grothflat$ of categories of the form $\Mod\left(\left(\BBB \otimes_\omega \lfp(\Qcoh(X))\right)^\op\right)$ where $\otimes_\omega$ denotes the tensor product of finitely cocomplete linear categories.
	\begin{proof}
		By \cite[Cor 6.9.12]{grothendieckdieudonne71} and \cite[Prop 75]{murfet06} we know that $\Qcoh(Y) \simeq \Ind(\lfp(\Qcoh(Y)))$. Then, as shown in the proof of \Cref{lem:reducingthebase}, we can write $\Qcoh(Y)^{\Qcoh(X)}$ as a pullback of exponentials of the form 
		\begin{equation}
			\begin{tikzcd}
				\Qcoh(Y)^{\Qcoh(X)} \arrow[rr] \arrow[d] &              &   \arrow[d] \Mod(\hat\BBB)^{\Qcoh(X)}     \\
				\Mod(\reallywidehat{\lfp(\Qcoh(Y))})^{\Qcoh(X)} \arrow[rr] &  			&\Mod( \hat\CCC )^{\Qcoh(X)}. 
			\end{tikzcd}
		\end{equation}
		In addition, as shown in the proof of \Cref{limitexp}, for any small linear category $\DDD$ we have that
		$$\Mod(\hat\DDD)^{\Qcoh(X)} \simeq (\Mod(k)[\obj])^{\Mod(\DDD^\op) \boxtimes \Qcoh(X)}.$$
		Now observe that, by the construction of $\otimes_\omega$ (see \cite[\S 2.4]{lopezfranco13}, \cite[\S5.1]{lowenramosgonzalezshoikhet18}), we have that
		\begin{equation*}
			\begin{aligned}
				\Mod(\DDD^\op) \boxtimes \Qcoh(X) &\simeq \Ind(\fpmod(\DDD^\op)) \boxtimes \Ind(\lfp(\Qcoh(X)))\\ 
				&\simeq \Ind(\fpmod(\DDD^\op) \otimes_\omega \lfp(\Qcoh(X))).
			\end{aligned}
		\end{equation*} 
		Then, by applying \Cref{locfinpresexponentiable}, we have that
		$$\Mod(\hat\DDD)^{\Qcoh(X)} \simeq \Mod((\fpmod(\DDD^\op) \otimes_\omega \lfp(\Qcoh(X)))^\op)$$
		as desired.
 
	\end{proof}	
\end{prop}

The presentation of $\Qcoh(Y)^{\Qcoh(X)}$ provided in \Cref{expgeneraldescription} involves the computation of the category $\fpmod(\lfp(\Qcoh(Y))^\op)$. We will show that we can avoid the computation of this category when $Y$ is noetherian. In order to do so, we need the following previous result. 
\begin{prop}\label{modulecategorybase}
	Let $\cg$ be an exponentiable Grothendieck category and $\AAA$ be a small linear category such that $\AAA^\op$ is ind-abelian in the sense of \cite[Def 1]{schappi14}, i.e. the ind-completion $\Ind(\AAA^\op)$ is an abelian category (and thus Grothendieck). Then
	\begin{equation}
		\Mod(\AAA)^\cg \simeq \Mod(k)[\obj]^{\Ind(\AAA^{\op})\boxtimes \cg}.
	\end{equation}
	\begin{proof}
		We have that
		\begin{equation}
			\begin{aligned}
				\grothflat(\cf \boxtimes \cg, \Mod(\AAA)) &\simeq \grothflat ( \cf \boxtimes \cg, \Mod(k)[\obj]^{\Ind(\AAA^{\op})})\\
				&\simeq \grothflat ((\cf \boxtimes \cg) \boxtimes \Ind(\AAA^{\op}), \Mod(k)[\obj])\\
				&\simeq \grothflat(\cf, \Mod(k)[\obj]^{\Ind(\AAA^{\op} )\boxtimes \cg} ),
			\end{aligned}
		\end{equation}
		where the first equality follows from \Cref{locfinpresexponentiable} and the third from the associativity and simetry of the tensor product of Grothendieck categories. 
	\end{proof}
\end{prop}

By means of \Cref{modulecategorybase}, we simplify the computation of the exponential $\Qcoh(Y)^{\Qcoh(X)}$ when $Y$ is noetherian.
\begin{prop}\label{expmodulecats}
	Let $X$ be a quasi-compact quasi-separated scheme and $Y$ a noetherian scheme. Then we have that $\Qcoh(X)^{\Qcoh(Y)}$ is the pullback in $\grothflat$ of a diagram of the form
	\begin{equation*}\label{eq:pullbacknoetherian}
		\begin{tikzcd}
			\Qcoh(Y)^{\Qcoh(X)} \arrow[rr] \arrow[d] &              &   \arrow[d] \Mod((\fpmod(\CCC^\op) \otimes_\omega \lfp(\Qcoh(X)))^\op)    \\
			\Mod\left(\left(\coh(Y)^\op \otimes_\omega \lfp(\Qcoh(X))\right)^\op\right) \arrow[rr] &  			&\Mod((\fpmod(\BBB^\op) \otimes_\omega \lfp(\Qcoh(X)))^\op). 
		\end{tikzcd}
	\end{equation*}
	In particular, if $X$ is also noetherian, we have that the left lower corner of the diagram can be written as
	$$\Mod((\coh(Y)^\op \otimes_{\mathsf{D}} \coh(X))^\op)$$
	where $\otimes_{\mathsf{D}}$ denotes the Deligne tensor product of abelian categories \cite{deligne90}.	
	\begin{proof}
		Observe that, as $Y$ is noetherian, $\Qcoh(Y) \simeq \Ind(\coh(Y)) \simeq \Mod(\coh(Y))[S]$, where $S$ is the family of morphisms $f:F \ra G$ in $\Mod(\coh(Y))$ such that $\Kern(f)$ and $\Cokern(f)$ are effaceable (see \cite[Thm 2.3]{krause15}). In particular, we can express $\Qcoh(Y)$ as a pullback diagram
		\begin{equation*}
			\begin{tikzcd}
				\Qcoh(Y) \arrow[rr] \arrow[d] &              &   \arrow[d] \Mod(\hat\CCC)   \\
				\Mod(\coh(Y))  \arrow[rr] &  			&\Mod(\hat\BBB). 
			\end{tikzcd}
		\end{equation*}
		as in the proof of \Cref{lem:reducingthebase}, but in this case without taking the free completion under finite limits of $\coh(Y)$. Then, we have that $\Qcoh(Y)^{\Qcoh(X)}$ is the pullback in $\grothflat$ of the diagram
		\begin{equation*}
			\begin{tikzcd}
				\Qcoh(Y)^{\Qcoh(X)} \arrow[rr] \arrow[d] &              &   \arrow[d] \Mod((\fpmod(\CCC^\op) \otimes_\omega \lfp(\Qcoh(X)))^\op)    \\
				\Mod(\coh(Y))^{\Qcoh(X)} \arrow[rr] &  			&\Mod((\fpmod(\BBB^\op) \otimes_\omega \lfp(\Qcoh(X)))^\op), 
			\end{tikzcd}
		\end{equation*}
		where the computation of the two elements on the right has been obtained as in the proof of \Cref{expgeneraldescription}.
		Because $\coh(Y)$ is abelian, so is $\coh(Y)^{\op}$. In particular, $\coh(Y)^{\op}$ is ind-abelian. Then, by \Cref{modulecategorybase}, we have that 
		\begin{equation*}
			\Mod(\coh(Y))^{\Qcoh(X)} \simeq \Mod(k)[\obj]^{\Ind(\coh(Y)^{\op}) \boxtimes \Qcoh(X)}.
		\end{equation*} 
	Moreover, we have that
	\begin{equation*}
	\begin{aligned}	
		\Ind(\coh(Y)^{\op}) \boxtimes \Qcoh(X) &\simeq \Ind(\coh(Y)^{\op}) \boxtimes \Ind(\lfp(\Qcoh(X)))\\ &\simeq \Ind\left(\coh(Y)^{\op} \otimes_{\omega} \lfp(\Qcoh(X)) \right).
	\end{aligned}
	\end{equation*}
		Therefore, if we now apply \Cref{locfinpresexponentiable}, we obtain
		\begin{equation*}
			\Mod(k)[\obj]^{\Ind(\coh(Y)^{\op}) \boxtimes \Qcoh(X)} \simeq \Mod\left(\left(\coh(Y)^{\op} \otimes_{\omega} \lfp(\Qcoh(X))\right)^\op\right),
		\end{equation*}
		as we wanted to show.
		
		Observe that, if $X$ is also noetherian, we have that $\lfp(\Qcoh(X)) \simeq \coh(X)$, which is an abelian category. As $\coh(Y)^\op$ is also an abelian category, the last claim of the statement follows from the fact that   
		\begin{equation*}
			\coh(Y)^\op \otimes_\omega \coh(X) \simeq \coh(Y)^\op \otimes_{\mathsf{D}} \coh(X) 
		\end{equation*}
	 	as shown in \cite[Thm 18]{lopezfranco13}.
	\end{proof}
\end{prop}

The following result will allow us to use \Cref{modulecategorybase} to compute $\Mod(A)^{\Qcoh(X)}$ where $A$ is an absolutely flat commutative algebra. 
\begin{prop}\label{fpabsolutelyflat}
	Let $A$ be an absolutely flat (also known as von Neumann regular) commutative algebra. Then $$\Mod(A) \simeq \Mod(\fpmod(A)).$$
	\begin{proof}
		From \cite[Thm 5.27]{kelly82basic}, we have that $\Mod(A) \simeq \Mod(\bar{A})$, where $\bar{A}$ denotes the linear Cauchy completion of the 1-object linear category defined by $A$. It hence suffices to show that $\bar{A} \simeq \fpmod(A)$. From \cite[Prop 3.4]{borceuxquinteiro96accessible}, we have that $\bar{A}$ can be recovered as $(\mathsf{Flat}(A^\op))_\omega$, the category of finitely presented objects of the category of flat $A^\op$-modules. In addition, as $A$ is commutative, we have that $A \simeq A^\op$. Consequently, $\bar{A} \simeq (\mathsf{Flat}(A))_\omega$. Furthermore, as $A$ is absolutely flat, we have that $\mathsf{Flat}(A) \simeq \Mod(A)$, and hence, $\bar{A} \simeq \fpmod(A)$, as we wanted to show. 
	\end{proof}
\end{prop}
\begin{exa}\label{exabsolutelyflat}
	Let $X$ be a noetherian scheme and $A$ an absolutely flat commutative algebra. Then, as a direct consequence of \Cref{fpabsolutelyflat} and \Cref{modulecategorybase} we have that
	\begin{equation}
		\Mod(A)^{\Qcoh(X)} \simeq \Mod\left( \left( \coh(X) \otimes_{\mathsf{D}} \fpmod(A)^{\op} \right) ^\op\right). 
	\end{equation}
\end{exa}

\section*{Acknowledgements}
	The first named author was supported by the Grant Agency of the Czech Republic project EXPRO 20-31529X and RVO: 67985840. The second named author is a postdoctoral researcher of the Fonds de la Recherche Scientifique - FNRS under Grant No. 32709538. She was a postdoctoral fellow of the Research Foundation - Flanders (FWO) under Grant No. 12T2619N during the time the authors were working on this article.
	
	The authors would like to thank an anonymous referee for the careful reading of the manuscript. They are also grateful to Wendy Lowen for the interesting discussions and to Andrea Gagna for his comments on the first draft of the introduction that led to an improvement of the references. The second named author would also like to thank Pieter Belmans for an interesting discussion about Hom-schemes. The first ideas that eventually led to this paper emerged during a research stay of the first named author at the University of Antwerp (Antwerp, Belgium). He would like to thank the University of Antwerp for the financial support that made his visit possible. An important part of the work presented in this article was done during a research stay of the second named author at Masaryk University (Brno, Czech Republic). She would like to thank Masaryk University and the Research Foundation - Flanders (FWO) for the financial support that made this research stay possible.

\bibliography{thebib}
\bibliographystyle{alpha}
\end{document}